\documentclass{article}
\usepackage{latexstyle}
\usepackage{multirow}
\usepackage{enumitem}
\usepackage{hyperref}
\usepackage{float}
\title{A practical guide to the recovery of wavelet coefficients from Fourier measurements}
\author{Milana Gataric\thanks{CCA, Centre for Mathematical Sciences, University of Cambridge, UK (m.gataric@maths.cam.ac.uk)} \and Clarice Poon\thanks{CCA, Centre for Mathematical Sciences, University of Cambridge, UK (cmhsp2@cam.ac.uk)}}

\date{April 2015; Revised March 2016}
\begin{document}
\maketitle

\begin{abstract}
In a series of recent papers (Adcock, Hansen and Poon, 2013, Appl.~Comput.~Harm.~Anal. 45(5):3132--3167), (Adcock, Gataric and Hansen, 2014, SIAM J.~Imaging~Sci.~7(3):1690--1723) and (Adcock, Hansen, Kutyniok and Ma, 2015, SIAM J.~Math.~Anal.~47(2):1196--1233), it was shown that one can optimally recover the wavelet coefficients of an unknown compactly supported function from pointwise evaluations of its Fourier transform via the method of generalized sampling. While these papers focused on the optimality of generalized sampling in terms of its stability and error bounds, the current paper explains how this optimal method can be implemented to yield a computationally efficient algorithm. In particular, we show that generalized sampling has a computational complexity of $\ord{M(N)\log N}$ when recovering the first $N$ boundary-corrected wavelet coefficients of an unknown compactly supported function from $M(N)$ Fourier samples. 
Therefore, due to the linear correspondences between the number of samples $M$ and number of coefficients $N$ shown previously, generalized sampling offers a computationally optimal way of recovering wavelet coefficients from Fourier data.
\end{abstract}

\section{Introduction}
The motivation behind the papers \cite{OptimalWaveletRecons,BAMGACHNonuniform1D,BAACHGKJM2Dwavelets} is that, there are countless applications in which one is concerned with the recovery of a function from finitely many pointwise evaluations of its Fourier transform. To name a few examples, these applications range from medical imaging \cite{guerquin2011fast,pruessmann1999sense,FesslerFastIterativeMRI} to electron microscopy \cite{Lawrence2012}, helium atom scattering \cite{Jardine_heliumscattering,jones_heliumscattering}, reflection seismology \cite{bleistein2000} and radar imaging \cite{Radar2005}. More precisely, one is tasked with the recovery of an unknown function $f\in \cL^2(\bbR^d)$ supported in a compact domain, from given samples of the form $\{\hat f(\omega): \omega\in \Omega\}$, where $\Omega$ is a countable set in $\bbR^d$ and the Fourier transform is defined as
$$
\hat f(\xi) = \int_{\bbR^d} f(x) \E^{2\pi i \xi\cdot x} \mathrm{d}x.
$$
The set of sampling points $\Omega$  either can be taken on an equidistant grid, which is known as \textit{uniform sampling}, or  can be an arbitrary countable set of certain density, which is known as \textit{nonuniform sampling}.  The listed papers described how, using generalized sampling, one can in principle optimally recover the (continuous) wavelet coefficients of a function, given its Fourier coefficients, and thereby directly exploit the superior approximation properties of wavelets. The present paper addresses the computational aspects of this generalized sampling method for both uniform and nonuniform sampling scenarios.

\subsection{The need for generalized sampling}

The classical approach to this recovery problem is to directly invert the Fourier sampling operator. In the case where $\Omega$ consists of equispaced points, one can simply apply the Nyquist-Shannon sampling theorem to compute a discrete inverse Fourier transform on the given samples, while in the case where $\Omega$ consists of nonequispaced points, one can approximate $f$ by inverting the frame operator associated with the Fourier frame $\left\{e^{2\pi i \omega \cdot}:\omega\in\Omega\right\}$ \cite{AldroubiGrochenigSIREV,BenedettoFrames}.
Another approach to nonuniform sampling is so-called gridding \cite{JacksonEtAlGridding}, a heuristic method which simply discretizes the inverse Fourier transform on a nonequidistant mesh by using certain density compensation factors.

In practice, it is always the case that only a finite number of samples is available, i.e. $\text{card}(\Omega)=M<\infty$, and hence, unless the original function is periodic and continuous, these direct Fourier approaches suffer from a number of drawbacks, including  slow convergence rates and artifacts such as the Gibbs phenomenon.  Therefore, these classical approaches are considered unsatisfactory \cite{unser2000sampling,GelbSong2014}.

In order to obtain a better reconstruction, one might seek an approximation represented in a different, more favorable system than the one used for sampling. In fact, a typical magnetic resonance (MR) image, for example, is known to be better represented by wavelets than by complex exponentials \cite{unser2000sampling}. The idea of basis change is not new in signal processing and it is typically referred to as \textit{generalized sampling}. It can be tracked down to the work of Unser and Aldroubi \cite{unser1994general,UnserHirabayashiConsist}, which was later extended by Eldar \cite{eldar2003sampling,eldar2005general}. Note that some of the early works of generalized sampling are sometimes referred to as consistent sampling. The idea of generalized sampling has also been applied in the work of Pruessmann et al. \cite{pruessmann1999sense} and Sutton, Noll and Fessler \cite{FesslerFastIterativeMRI} for the recovery of voxel coefficients from Fourier samples and also in the work of Hrycak and Gr\"{o}chenig \cite{hrycakIPRM} for the recovery of polynomial coefficients from Fourier samples. Furthermore, the efficient recovery of trigonometric polynomial coefficients of bandlimited functions from its nonuniform samples, or equivalently, the recovery of Dirac coefficients of compactly supported $\mathcal{L}^2$ functions from nonuniform Fourier samples, was considered in works by Feichtinger, Gr\"ochenig and Strohmer \cite{Feichtinger_efficientnon}. 

These approaches were formalized under the framework of Adcock and Hansen \cite{BAACHShannon} in the setting of a Hilbert space $\cH$, where it was shown that one can recover an unknown element $f\in\cH$ in terms of any frame $\br{\varphi_j}_{j\in\bbN}$ when given samples of the form $(\ip{f}{\psi_j})_{j\in\bbN}$ where $\br{\psi_j}_{j\in\bbN}$ is an arbitrary frame in $\cH$. More specifically, given the first $M$ samples  $\br{\ip{f}{\psi_j}: j=1,\ldots, M}$,  one can approximate the first $N$ reconstruction coefficients $\br{\ip{f}{\varphi_j}: j=1,\ldots, N}$ in a stable and quasi-optimal manner (see Section \ref{sec:background}  for the definition) provided that $M$ is sufficiently large with respect to $N$.

The setting  where $\br{\varphi_j}_{j\in\bbN}$ forms a wavelet basis and $\br{\psi_j}_{j\in\bbN}$ forms a Fourier basis was considered in \cite{OptimalWaveletRecons} and \cite{BAACHGKJM2Dwavelets}  in one and two dimensions respectively. It was shown that  a linear scaling between the number of Fourier samples  $M$ and the number of  wavelet coefficients $N$ is sufficient for stable and quasi-optimal recovery. In other words, it suffices to take  $M=\ord{N}$ Fourier samples of a function in order to stably quasi-optimally recover its first $N$ wavelet coefficients. A more general setting where  $\br{\psi_j}_{j\in\bbN}$ forms a (weighted) Fourier frame was considered in \cite{BAMGACHNonuniform1D}, in one dimension, where it was shown that it is sufficient to take Fourier samples from the interval $[-K,K]$ with $K=\ord{N}$ in order to recover the first $N$ wavelet coefficients. Therefore, as long as the number of samples $M$ scales linearly with the sampling bandwidth $K$, the relationship $M=\ord{N}$ is also sufficient for stable and quasi-optimal recovery in this more general setting.

\subsection{Contribution of this paper}
The method of generalized sampling can be recast as the linear system
\be{\label{GSmatrix}
\begin{pmatrix}
\ip{\varphi_1}{\psi_1} & \ip{\varphi_2}{\psi_1} & \ldots & \ip{\varphi_N}{\psi_1} \\ 
\ip{\varphi_1}{\psi_2} & \ip{\varphi_2}{\psi_2} & \ldots & \ip{\varphi_N}{\psi_2} \\ 
\vdots & \vdots & \ddots & \vdots \\ 
\ip{\varphi_1}{\psi_M} & \ip{\varphi_2}{\psi_M} & \ldots & \ip{\varphi_N}{\psi_M}
\end{pmatrix}
\begin{pmatrix}
 \alpha_1\\
 \alpha_2\\
 \vdots\\
 \alpha_N 
\end{pmatrix}
=
\begin{pmatrix}
\ip{f}{\psi_1}\\
\ip{f}{\psi_2}\\
\vdots\\
\ip{f}{\psi_M}
\end{pmatrix},
}
which is to be solved for the least-squares solution $\{\alpha_j: j=1,\ldots, N\}$ that approximates the desired coefficients $\br{\ip{f}{\varphi_j}: j=1,\ldots, N}$, given the $M$ samples $\br{\ip{f}{\psi_j}: j=1,\ldots, M}$. In the general case, solving this system has a computational complexity of $\ord{NM}$. In this paper we show that in the case of recovering  wavelet coefficients from  Fourier samples, the computational  complexity of generalized sampling is only $\ord{M\log N}$, since in this case the cost of applying the matrix from \R{GSmatrix} and its adjoint is only $\ord{M\log N}$. In fact, in uniform sampling as well as in nonuniform sampling where $M=\ord{K}$, the computational complexity is simply $\ord{N\log N}$, due to the linear correspondences derived in the aforementioned papers.   Therefore, one can attain the superior reconstruction qualities via generalized sampling at a computational cost that is on the same order as the computational cost of the classical methods that are based on a simple discretization of the inverse Fourier transform. 

This paper will describe the computational issues relating to the recovery of boundary-corrected wavelet coefficients \cite{cohen1993wavelets} when given either uniform or nonuniform Fourier samples. Although we shall only address the reconstruction of coefficients in one and two dimensions, the techniques described in this paper can readily be applied to higher dimensional cases. Generalized sampling may also be efficiently implemented with wavelets satisfying other boundary conditions (such as periodic or symmetric boundary conditions), however, here we consider the boundary-corrected wavelets of \cite{cohen1993wavelets} since such wavelets preserve vanishing moments at the domain boundaries and form unconditional bases on function spaces of certain regularity on bounded domains.

While in this paper we mainly focus on the linear recovery model \R{GSmatrix}, the computational aspects analyzed in this paper arise in various other nonlinear recovery schemes such as the $\ell_1$-minimization schemes introduced in \cite{gs_l1,BAACHGSCS,adcockbreaking}. Namely, whenever one wants to recover wavelet coefficients from Fourier measurements, one needs fast computations involving the same matrix as the one appearing in \R{GSmatrix}, as well as the fast computations involving its adjoint. Hence, the algorithms described in this paper can readily be applied to yield efficient implementations of these other  nonlinear recovery schemes.

We remark that although the efficient recovery of wavelet coefficients from Fourier samples was already considered in \cite{guerquin2011fast}, that work did not provide a comprehensive explanation of the implementation using boundary-corrected wavelets and various types of Fourier sampling. Moreover, the purpose of this work is to provide a practical guide to the use of generalized sampling with wavelets and Fourier sampling. 
A MATLAB implementation of our algorithm is available at  \url{http://www.damtp.cam.ac.uk/research/afha/code/}.

\subsection{Outline}
Section \ref{sec:background} describes the framework of generalized sampling. Section \ref{s:wavelets_from_fourier}  describes the special case of generalized sampling which produces wavelet reconstructions from Fourier samples. The wavelet reconstruction space and Fourier sampling space are defined in Section \ref{ss:rs} and Section \ref{ss:ss}, respectively, and additionally, in Section \ref{ss:ssr}, we recall the results providing the linear stable sampling rate for this particular pair of spaces. In Section \ref{s:algorithm_1d}, we deliver the detailed reconstruction algorithm for the one-dimensional case together with the analysis of its computational complexity. The two-dimensional case is presented in Section \ref{s:algorithm_2d} in a general setting, and the simplifications arising from uniform sampling are presented in Section \ref{s:algorithm_2d_simp}. In Section \ref{s:num_ex}, we use the presented algorithm in a few  numerical examples.

\subsection{Notation}
Let $\cH$ denote a Hilbert space with norm $\norm{\cdot}$ and inner product $\ip{\cdot}{\cdot}$. A set $\{\phi_m\}_{m\in\bbN} \subseteq \cH$ forms a frame for $\cH$ if there exist constants $A,B>0$ such that 
\be{\label{frame_ineq}
A \|f\|^2 \leq \sum_{m\in\bbN} \abs{\ip{f}{\phi_m}}^2 \leq B \|f\|^2, \quad \forall f \in \cH.
}
For any subspace $\cT\subseteq\cH$, let $Q_{\cT}$ denote the orthogonal projection onto $\cT$.
Throughout the paper, we shall consider the Hilbert space of square-integrable $\cL^2$ functions which are compactly supported on the cube $[0,1]^d\subseteq\bbR^d$ in the physical domain, for dimension $d\in\bbN$.
For a point in the frequency domain $\omega\in\bbR^d$ let 
$$
e_{\omega}(x)=\E^{2\pi\I\omega \cdot x}\chi_{[0,1]^d}(x),\quad x\in\bbR^d.
$$
Hence, for a function $f\in\cL^{2}([0,1]^d)$, we can write
$$\hat{f}(\omega)=\ip{f}{e_{\omega}}$$
for the Fourier transform of $f$ at $\omega$. In particular, if \R{frame_ineq} holds with $\phi_m=e_{\omega_m}$, we say that $\{e_{\omega_m}\}_{m\in\bbN}$ is a Fourier frame for $\cL^{2}([0,1]^d)$.
For a matrix $A\in \bbC^{M\times N}$, we use the notation $A=(A_{m,n})_{n=1,m=1}^{N,M}$ and denote its pseudoinverse by $A^{\dagger}$ and its adjoint by $A^*$.

\section{Generalized sampling}\label{sec:background}  
This section recaps the main features of the generalized sampling method described in \cite{BAACHShannon,BAACHOptimality}.
Let $\cR,\cS$ be closed subspaces of $\cH$ such that $\cR\cap \cS^\perp=\{0\}$ and $\cR + \cS^\perp$ is a closed subspace of $\cH$. Let $\br{\varphi_n}_{n\in\bbN}$ and $\br{\psi_m}_{m\in\bbN}$ be frames for $\cR$ and $\cS$ respectively.
For $N,M\in\bbN$, let 
$$
\cR_N = \mathrm{span}\br{\varphi_n:n=1,\ldots, N}, \quad \cS_M = \mathrm{span}\br{\psi_m:m=1,\ldots, M}
$$
denote the finite dimensional reconstruction and the sampling space respectively. For an unknown function $f\in\cH$, we seek an approximation $\tilde f \in \cR_N$ to $f$ using only the finite set of samples $\left\{\ip{f}{\psi_m}:m=1,\ldots, M\right\}$. 

It was shown in \cite{BAACHOptimality} that for each $N\in\bbN$, there exists $m_0\in\bbN$ such that for all $M\geq m_0$, there exists a unique element $\tilde f\in\cR_N$ satisfying
\be{\label{eq:gs}
\ip{S_M \tilde f}{\varphi_n} = \ip{S_M f}{\varphi_n}, \quad n=1,\ldots, N.
}
where $S_M:\rH\rightarrow\cS_M$, $g\mapsto\sum_{m=1}^{M}\ip{g}{\psi_m}\psi_m$. The approximation $\tilde f$ is referred to as the generalized sampling reconstruction. Furthermore, this reconstruction satisfies the sharp bounds
\be{\label{bounds}
\|{\tilde f - f}\| \leq C_{N,M}\nm{Q_{\cR_N} f - f}, \qquad \nmu{\tilde f} \leq C_{N,M} \nm{f},
}
where
$$
C_{N,M}^{-1} = \inf_{\substack{u\in\cR_N\\ \nm{u}=1}} \nm{Q_{\cS_M} u} >0.
$$
Note that the approximation $\tilde f$ achieves, up to a constant of $C_{N,M}$, the same error bound as $Q_{\cR_N}f$, which is the orthogonal projection onto $\cR_N$ and therefore the best possible approximation to $f$ in $\cR_N$. Hence, $\tilde f$ may be said to be \textit{quasi-optimal}.

Due to \R{eq:gs}, by defining the generalized sampling operator
\be{\label{eq:gs_op}
G^{[N,M]}: \bbC^N \to \bbC^M, \qquad \alpha \mapsto \left(\ip{\sum_{n=1}^N \alpha_n \varphi_n}{\psi_m}\right)_{m=1}^M,
}
the solution $\tilde f$ can be written as 
$$\tilde f = \sum_{j=1}^N \alpha_j \varphi_j,$$
where $\alpha=(\alpha_j)_{j=1}^N\in  \bbC^N$ is the least  squares solution to
\be{\label{eq:least_sq}
G^{[N,M]}(\alpha) = (\ip{f}{\psi_j})_{j=1}^M.
}
Note that this is the same as the system written in \R{GSmatrix}.
Furthermore, the condition number of the operator $G^{[N,M]}$, which indicates reconstruction stability, can be shown to be
$$
D_{N,M}= \left( \inf_{u\in\cR_N} \sum_{m=1}^M \abs{\ip{u}{\psi_m}}^2\right)^{-1/2}>0.
$$

The accuracy and stability of generalized sampling are therefore governed by $C_{N,M}$ and $D_{N,M}$ and it is imperative to determine the correct scaling between the number of reconstruction vectors $N$ and the number of samples $M$ to ensure boundedness of these two values. 
For a given number of reconstruction vectors $N$ and some $\theta \in (0,\infty)$, the stable sampling rate is defined as a minimal number of sampling vectors $M$ providing a stable and quasi-optimal reconstruction. Formally, the stable sampling rate is given as 
\be{\label{ss_rate}
\Theta(N,\theta) = \min\br{M\in\bbN : \max\br{C_{N,M}, D_{N,M}} \leq \theta}.
}
The stable sampling rate for particular choices of wavelet reconstruction spaces and Fourier sampling spaces is discussed further in Section \ref{ss:ssr}.

\subsection{The computational challenges}\label{ss:challenges}
There are two aspects to the  implementation of generalized sampling in (\ref{eq:gs}):
\begin{enumerate}
\item We need to understand the stable sampling rate $M = \Theta(N,\theta)$ required for quasi-optimal and stable reconstructions.  
\item Having determined the appropriate finite section matrix $G^{[N,M]}$,  it remains to solve least-squares problem \R{eq:least_sq}. Note that given any matrix $G\in\bbC^{M\times N}$, there are many iterative schemes for computing the least-squares solution $G^\dagger g$. For example, one can apply the conjugate gradient method to the corresponding normal equations \cite{trefethen1997numerical}.  The efficiency of these iterative methods is always dependent on the computational complexity of applying $G$ and $G^*$, and in the worst case, this is $\ord{NM}$.
\end{enumerate} 
In the case of recovering wavelet coefficients from Fourier samples, the first challenge detailed above has been resolved in 
 \cite{OptimalWaveletRecons} where the sampling vectors form a basis, and in \cite{BAMGACHNonuniform1D}, where the sampling vectors form a 
(weighted) Fourier frame. The main results from these works will be summarised in  Section \ref{ss:ssr}. The purpose of this paper is to show that the computational complexity of applying $\rG^{[N,M]}$ and its adjoint $(\rG^{[N,M]})^*$ is $\ord{M\log N}$. Thus, because of the linear correspondence between the number of samples and reconstruction coefficients revealed in  \cite{OptimalWaveletRecons} and \cite{BAMGACHNonuniform1D}, generalized sampling offers a computationally optimal way of recovering wavelet coefficients from Fourier data.

\rem{[Generalized sampling and $\ell^1$ minimization] Based on ideas of compressed sensing introduced in \cite{donohoCS,candes2006robust},
to recover $f\in \cL^2$ from a reduced number of noisy Fourier measurements $y = \{\hat f(\omega_j): j\in\Omega\}+\eta$, where $\nm{\eta}_{\ell^2}\leq \delta$ and $\Omega$ is some finite index set, it was proposed in \cite{BAACHGSCS,adcockbreaking} to solve
\be{\label{eq:inf_dim_cs}
x \in \mathrm{arg} \inf_{z\in\ell^1(\bbN)} \nm{z}_{\ell^1} \text{ subject to } \nm{ G_\Omega z - y}_{\ell^1}\leq \delta,
}
where  $$G_\Omega = (\ip{\varphi_j}{e^{2\pi i k \cdot}})_{j\in\bbN,k\in\Omega}$$ is a semi-infinite dimensional generalized sampling matrix.  Solutions to this minimization problem were analysed in detail in \cite{adcockbreaking,gs_l1}, and we refer to those works for theoretical error estimates. However, we simply mention that if $\Omega = \br{j : j=-M,\ldots, M}$ and $\delta=0$, it was proved in \cite{gs_l1} that if $\br{\varphi_j}_{j\in\bbN}$ is a wavelet basis with sufficiently many vanishing moments, then  any solution $x$ to (\ref{eq:inf_dim_cs}) will satisfy 
$$
\nm{f - \sum_{j=1}^\infty x_j \varphi_j}_{L^2} =\ord{  \sum_{j=N'}^\infty \abs{\ip{f}{\varphi_j}}}
$$
where $N' = cN$ for some constant $c\in (0,1]$ which is independent of $N$. Thus, the reconstruction error is governed by the decay in the wavelet coefficients of $f$, rather than the decay in its Fourier coefficients.
Computationally, one can show that any solution to 
 \be{\label{eq:gsl1_fin}
\inf_{\alpha\in\bbC^K}\nm{\alpha}_{\ell^1} \text{ subject to } \nm{ G^{[K,M]}(\alpha) - \beta }_{\ell^2}\leq \delta
}
converges in an $\ell^1$ sense to solutions of (\ref{eq:inf_dim_cs}) as $K\to \infty$  (see \cite{BAACHGSCS}) and for the Fourier-wavelets case, it suffices to let $K=\ord{M}$ \cite{adcockbreaking}. There are again numerous algorithms such as \cite{BergFriedlander:2008, spgl1:2007} for solving problems of the form (\ref{eq:gsl1_fin}) and, similarly as for \R{eq:gs}, the efficiency of these algorithms  amounts to the efficiency of applying $G^{[K,M]}$ and its adjoint. 
}

\section{The wavelet reconstruction from the Fourier samples}\label{s:wavelets_from_fourier}

As mentioned before, in this paper we consider the particular reconstruction space consisting of $N$ boundary-corrected wavelets and the sampling space consisting of $M$ Fourier-type exponentials. In this section, we define the corresponding spaces and give the stable sampling rate for this pair.

\subsection{The reconstruction space}\label{ss:rs}

We shall consider the Daubechies wavelets of $a$ vanishing moments, on the interval $[0,1]$, with the boundary wavelet construction introduced in \cite{cohen1993wavelets}. Following convention, we assume that the scaling function $\phi$ and the wavelet $\psi$ are supported on $[-a+1,a]$ and denote boundary scaling functions at the endpoints $0$ and $1$ by 
$$\phi^0_k,\ \phi^1_k, \qquad k=0,\ldots,a-1$$
and boundary wavelets at the endpoints $0$ and $1$ by
$$\psi^0_k,\ \psi^1_k, \qquad k=0,\ldots,a-1.$$
We denote the boundary-corrected scaling functions on the interval $[0,1]$ by
\bes{
\phi^{[0,1]}_{j,k}(x) = \left \{ \begin{array}{ll} 2^{j/2} \phi (2^j x - k) & a \leq k < 2^j - a 
\\ 
2^{j/2} \phi^{0}_{k} (2^j x) & 0 \leq k < a 
\\
2^{j/2} \phi^{1}_{2^j-k-1}(2^jx) & 2^j - a  \leq k < 2^j,
\end{array} \right .
}
and similarly for the wavelet functions $\psi^{[0,1]}_{j,k}$.

\paragraph{One dimension.} Here we define the reconstruction space for dimension $d=1$, consisting of boundary-corrected wavelets on the interval $[0,1]$ of $a$ vanishing moments.  Let $J$ be such that $J \geq \log_2(2a)$. The set
\be{\label{eq:daub_boundary}
\br{\phi^{[0,1]}_{j,k}: k=0,\ldots, 2^j-1} \cup \left(\bigcup_{j\geq J}\br{\psi^{[0,1]}_{j,k}: k=0,\ldots, 2^j-1} \right)
}
forms an orthonormal basis for $\cL^2([0,1])$. Furthermore, by defining
\bes{
W_j^{[0,1]} = \mathrm{span}\br{\psi^{[0,1]}_{j,k}: k=0,\ldots, 2^j-1},\ V_j^{[0,1]}= \mathrm{span}\br{\phi^{[0,1]}_{j,k}: k=0,\ldots, 2^j-1},
}
we have that
$\cL^2([0,1]) = V_J^{[0,1]} \oplus \left(\mathop{\oplus}_{J\leq j} W_j^{[0,1]} \right)$. We now order elements of (\ref{eq:daub_boundary}) in increasing order of wavelet scale, and we denote by $\cW_N$ the space spanned by the first $N$ wavelets. Thus, for fixed $R>J\geq \log_2(2a)$ and $N=2^R$ we have the reconstruction space
\be{
\label{eq:rec_space_wavelets}
\cW_{N} = V^{[0,1]}_{J} \oplus W^{[0,1]}_J \oplus \cdots \oplus W^{[0,1]}_{R-1}.
}
Note that, by  construction, $\text{dim}(\cW_{N})=N=2^R$ and also $\cW_{N}=V^{[0,1]}_{R}$. 

For a function $f\in\cW_{N}$ we can write
\bes{
f(x) = \sum_{k=0}^{2^J-1} c_{J,k} \phi^{[0,1]}_{J,k} (x) + \sum_{j=J}^{R-1} \sum_{k=0}^{2^j-1} d_{j,k} \psi^{[0,1]}_{j,k} (x)
}
and also
\bes{
f(x) = \sum_{k=0}^{2^R-1} c_{R,k} \phi^{[0,1]}_{R,k} (x) 
}
for some scaling coefficients $c_{j,k}$ and some detail coefficients $d_{j,k}$. We recall that, given the scaling coefficients $\{c_{R,k}:k=0,\ldots,2^R-1\}$, it is possible compute  the scaling coefficients $\{c_{J,k} : k=0,\ldots,2^J-1\}$ and detail coefficients $\{d_{j,k} : k=0,\ldots,2^j-1, \ j=J,\ldots,R-1  \}$, and vice versa. This can be done by the discrete boundary-corrected Forward Wavelet Transform (FWT), which we denote by $W$ and by $\textit{\textbf{W}}$ in two dimensions. The reverse operation is performed by the discrete boundary-corrected Inverse Wavelet Transform (IWT), denoted by $W^{-1}$, in one, and by $\textit{\textbf{W}}^{-1}$ in two dimensions.

\paragraph{Two dimensions.}
The boundary-corrected wavelet basis for $\cL^2([0,1]^2)$ is defined as follows. We introduce the following two-dimensional functions:
\eas{
\Phi_{j,(k_1,k_2)}(x_1,x_2) =  \phi^{[0,1]}_{j,k_1}(x_1) \phi^{[0,1]}_{j,k_2}(x_2), \quad \Psi_{j,(k_1,k_2)}^1(x_1,x_2) =  \phi^{[0,1]}_{j,k_1}(x_1) \psi^{[0,1]}_{j,k_2}(x_2), \\
\Psi_{j,(k_1,k_2)}^2(x_1,x_2) = \psi^{[0,1]}_{j,k_1}(x_1) \phi^{[0,1]}_{j,k_2}(x_2), \quad  \Psi_{j,(k_1,k_2)}^3(x_1,x_2) = \psi^{[0,1]}_{j,k_1}(x_1) \psi^{[0,1]}_{j,k_2}(x_2).
}
As before, take $J\geq \log_2(2a)$, and let
\bes{
\Upsilon_J^0=\br{\Phi_{J,(k_1,k_2)}: 0\leq k_1,k_2 \leq 2^J-1 }
}
and 
\bes{
\Upsilon^{i}_j = \br{\Psi^i_{j,(k_1,k_2)}: 0\leq k_1,k_2 \leq 2^j-1 }, \quad  \ j\in\bbN,\ j\geq J, \ i=1,2,3.
}
The set
\begin{equation}\label{eq:daub_boundary2d}
 \Upsilon_J^0  \cup \left(\bigcup_{j\geq J}\br{\Upsilon^{i}_j : i=1,2,3 }\right)
\end{equation}
forms an orthonormal basis for $L^2([0,1]^2)$.
We now order the basis elements of (\ref{eq:daub_boundary2d}) in increasing order of wavelet scales so that we can write
$$
\{\varphi_{n_1,n_2}\}_{n_1,n_2 \in\bbN^2} = \begin{pmatrix}
&\Upsilon^0_J &\Upsilon_J^1  & &\Upsilon_{J+1}^1 & \hdots\\
&\Upsilon_J^2  & \Upsilon_J^3 \\
\\
&&\!\!\!\!\!\!\!\! \Upsilon_{J+1}^2 & &\Upsilon_{J+1}^3\\
&&\!\!\!\!\!\!\!\! \vdots & & & \ddots
\end{pmatrix}.
$$
Let  $\cW_{N^2}$ denote the space spanned by the first $N\times N$ wavelets via this ordering, so that $\cW_{N^2}=\mathrm{span}\br{\varphi_{n_1,n_2}: 0\leq n_1,n_2 \leq N-1}$.
For $N=2^R$ and  $R>J\geq \log_2(2a)$, we have
\be{
\label{eq:rec_space_wavelets_2d}
\cW_{N^2} = \mathrm{span}\ \Upsilon_J^0 \oplus \left( \oplus_{i=1}^{3} \oplus_{j=J}^{R-1} \mathrm{span}\ \Upsilon_j^i \right ) = \mathrm{span}\ \Upsilon_R^0, 
}
which is the two-dimensional reconstruction space.

\subsection{The sampling space}\label{ss:ss}

Let $K>0$ be a sampling bandwidth and $Z_K\subseteq \bbR^d$ be a closed, simply connected set in the frequency domain such that $\max_{z\in Z_K}\abs{z}_{\infty}=K$. For finite $M=M(K)$, let 
$$
\Omega_M=\{\omega_{m} \ :\ m=1,\ldots,M \} \subseteq Z_K,
$$
be a given sampling scheme, a set of distinct frequency points in $Z_K$ which are permitted to be completely nonuniform.
The sampling space that we shall consider in this paper is given by
\be{\label{eq:wframe_sampling_space}
\cE_{M} = \mathrm{span} \br{ \sqrt{\mu_m} e_{\omega_m} :  \omega_m \in \Omega_M},
}
with the weights $\mu_m$  defined as measures of the Voronoi regions associated to the sampling points, i.e,
\be{\label{eq:weights}
\mu_{\omega_m} = \int_{Z_K} \chi_{V_{\omega_m}}(x)  \D x,
}
where
$$
V_{\omega_m} = \left\{ v\in Z_K : \forall n \neq m, \abs{v-\omega_m}_{\ell^1} \leq \abs{v-\omega_n}_{\ell^1} \right\}.
$$
The use of Voronoi weights is common in both practice and nonuniform sampling theory \cite{AldroubiGrochenigSIREV,spyralScienceDirect}.

As shown in \cite{AGH2DNUGS}, it is crucial that the sampling scheme $\Omega_M$ satisfy an appropriate density condition in the frequency region $Z_K$; i.e. it is crucial that
\be{\label{eq:density_con}
\delta_{Z_K}(\Omega_M) = \sup_{v\in Z_K} \inf_{\omega\in\Omega_M} \abs{v-\omega}_{\ell^1} <1/2,
}
but due to use of weights,  $\Omega_M$ is allowed to cluster arbitrarily as $K, M(K)\rightarrow\infty$. Namely, the density condition ensures that weighted exponentials in $\cE_{M}$ give a weighted Fourier frame for $\cL^2([0,1]^d)$ as $K, M(K)\rightarrow\infty$, i.e. \R{frame_ineq} is satisfied with $\phi_m= \sqrt{\mu_m} e_{\omega_m}$. This, on the other hand, ensures stable sampling.

The particular examples of nonuniform sampling schemes that we consider in the present paper are jittered and log sampling schemes in one dimension, and spiral and radial sampling schemes in two dimensions. The construction of these particular sampling schemes, such that they satisfy density condition \R{eq:density_con}, was described in \cite{BAMGACHNonuniform1D} for the one-dimensional case, and in \cite{AGH2DNUGS} for the two-dimensional case.

A special case of sampling that we also consider in this paper is uniform sampling. Here, the sampling scheme is taken on an equidistant grid with a fixed spacing. For a given spacing $\epsilon\in (0,1]$, the uniform sampling space in one dimension is defined as
\be{\label{eq:uniform_space_1d}
\cE_{M} = \mathrm{span}\br{\epsilon e_{\epsilon k} : k= -\lceil M/2\rceil, \ldots, \lceil M/2\rceil -1}
}
and in two dimensions as
\be{\label{sampling_space_uniform2d}
\cE_{M^2} = \mathrm{span}\br{\epsilon_1 \epsilon_2 e_{\epsilon_1 k_1} e_{\epsilon_2 k_2}: -\lceil M/2\rceil \leq k_1, k_2 \leq \lceil M/2\rceil -1},
}
for some $\epsilon_i\in (0,1]$, $i=1,2$.
It is known that exponentials in the uniform sampling space lead to a tight Fourier frame for $\cL^2([0,1]^d)$. As we shall see, for this special sampling scenario, the two-dimensional computations can be considerably simplified.

\rem{Since the clustering of the sampling points is allowed, the stable sampling rate in the case of nonuniform sampling is defined by
\be{\label{ss_rate_K}
\Theta(N,\theta) = \min\br{K>0 : \max\br{C_{N,M(K)}, D_{N,M(K)}} \leq \theta},
}
for a given number of reconstruction vectors $N\in\bbN$ and some $\theta \in (0,\infty)$.
Therefore, rather than a minimal number of sampling vectors $M$, as it is the case in \R{ss_rate}, here we seek for a minimal sampling bandwidth $K$ that provides a stable and accurate reconstruction.
}

\subsection{The stable sampling rate}\label{ss:ssr}

To emphasize the benefits of  using generalized sampling  to recover wavelet coefficients from Fourier samples, in this section, we recall the main results of \cite{OptimalWaveletRecons,BAMGACHNonuniform1D} and present a numerical example which demonstrates the superior convergence rates that generalized sampling offers.

To begin, generalized sampling is particularly effective in the regime of recovering wavelet coefficients from Fourier samples because a linear scaling between the number of Fourier samples $M$ (or sampling bandwidth $K$) and the number of wavelets $N$ is enough for stable and accurate recovery via generalized sampling. In one dimension we have the following.

\begin{theorem}[\cite{OptimalWaveletRecons,BAMGACHNonuniform1D}]\label{t:samp_rates} 
Let $\cW_N$ be the reconstruction space consisting of $N$ boundary-corrected Daubechies wavelets given by \R{eq:rec_space_wavelets}. 

\begin{itemize}
 \item[$(i)$] Let the sampling space $\cE_M$ be given by \R{eq:uniform_space_1d}.  Then for any $\gamma\in(0,1)$ there exists a constant $c_0=c_0(\gamma)$ such that
$$
\min\left\{M : \max\left\{C_{N,M},D_{N,M}\right\} \leq \gamma \right\} \leq c_0 N.
$$  
 \item[$(ii)$] More generally, let $\Omega_M$ be contained in the interval $[-K,K]$ and such that $\delta=\delta_{[-K,K]}(\Omega_M)<1/2$, i.e. let the sampling space $\cE_M$ be given by \R{eq:wframe_sampling_space}. Then for any $\gamma\in(0,1-\delta)$ there exists a constant $c_0=c_0(\gamma)$ such that
$$
\min\left\{K : \max\left\{C_{N,M},D_{N,M}\right\} \leq \frac{1+\delta}{1-\delta-\gamma}\right\} \leq c_0 N.
$$ 
\end{itemize}
\end{theorem}

This result characterizes the Fourier samplings which permit stable and accurate recovery in spaces of boundary-corrected wavelets. Both parts of this theorem provide estimates on a stable sampling rate of the type 
$$
\Theta(N,c_\gamma)=\cO(N),
$$
where  $\Theta$ is defined by \R{ss_rate} for uniform sampling, and  by \R{ss_rate_K} for nonuniform sampling.  Let us mention here that similar higher-dimensional results were obtained for the uniform Fourier samples in \cite{BAACHGKJM2Dwavelets}, while for nonuniform sampling, the higher-dimensional case remains an open problem.

The result given in \cite{mallat2008wavelet} states that if  $f\in \cH^s[0,1]$, $s\in(0,a)$, then  $ \norm{f-Q_{\cW_N} f} =  \mathcal{O}(N^{-s})$, where $\cH^s$ denotes the usual Sobolev space and $\cW_N$ is the space of $N$ boundary-corrected wavelets with $a$ vanishing moments. This result in combination with Theorem \ref{t:samp_rates} gives the following.

\begin{corollary}[\cite{OptimalWaveletRecons_sampta,BAMGACHNonuniform1D}]
Let $\mathcal{W}_N$ be the reconstruction space consisting of $N$ boundary-corrected Daubechies wavelets of $a$ vanishing moments  \R{eq:rec_space_wavelets}. Let $f\in \cH^s[0,1]$ with $s\in (0,a)$ be an arbitrary function.
\begin{itemize}
 \item[$(i)$] Let the sampling space $\cE_M$ be given by \R{eq:uniform_space_1d}, then the generalized sampling solution $\tilde f$ implemented with $M= \Theta(N, \gamma)$ samples satisfies
$$
\|f- \tilde f\| = \mathcal{O}(M^{-s}).
$$
 \item[$(ii)$] Let the sampling space $\cE_M$ be given by \R{eq:wframe_sampling_space}, then the generalized sampling solution $\tilde f$ implemented with the sampling bandwidth $K(M)= \Theta(N, (1+\delta)/(1-\delta-\gamma))$  satisfies
$$
\|f- \tilde f\| = \mathcal{O}(K^{-s}).
$$
\end{itemize}\label{cor:rates}
\end{corollary}

Hence, up to constant factors, generalized sampling obtains optimal convergence rates in terms of the number of sampling points $M$ (or the sampling bandwidth $K$), when reconstructing smooth functions with boundary-corrected Daubechies wavelets. 

\begin{figure}[h]
\begin{center}
\begin{tabular}{cccc}
\multicolumn{1}{l}{\hspace{1.3cm} Equispaced sampling \hspace{1.8cm} Jittered sampling \hspace{2.3cm} Log sampling} \vspace{0.1cm} \\
 \!\!\!\!{\rotatebox{90}{\hspace{2.4cm}$-\log\|f- \tilde f\|$}} \! \includegraphics[scale=0.67]{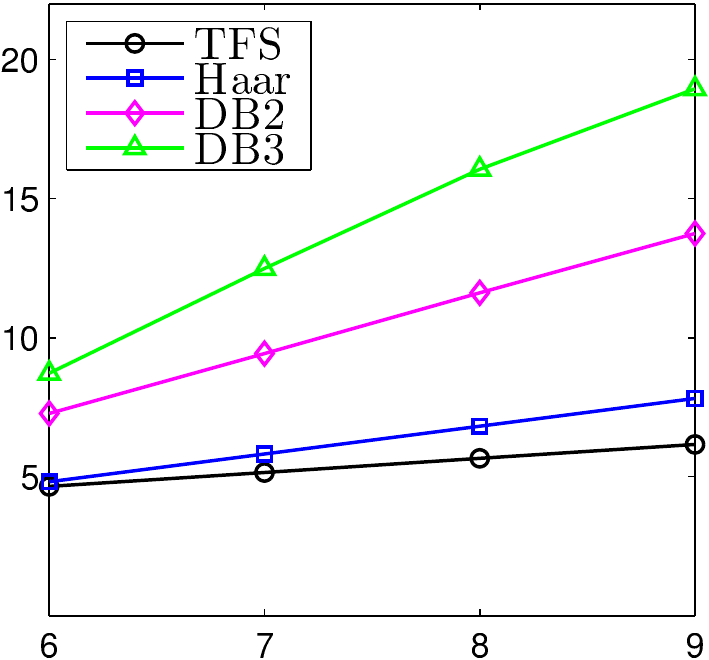}   \includegraphics[scale=0.67]{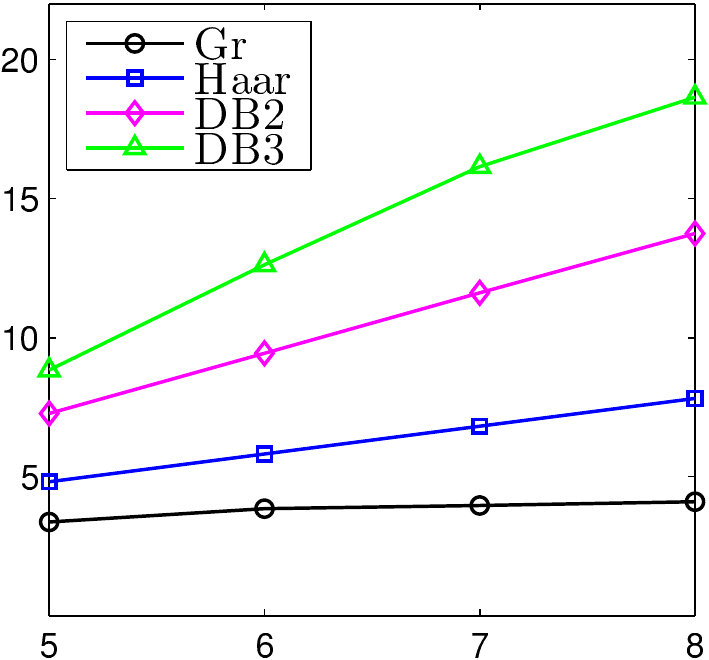}  \includegraphics[scale=0.67]{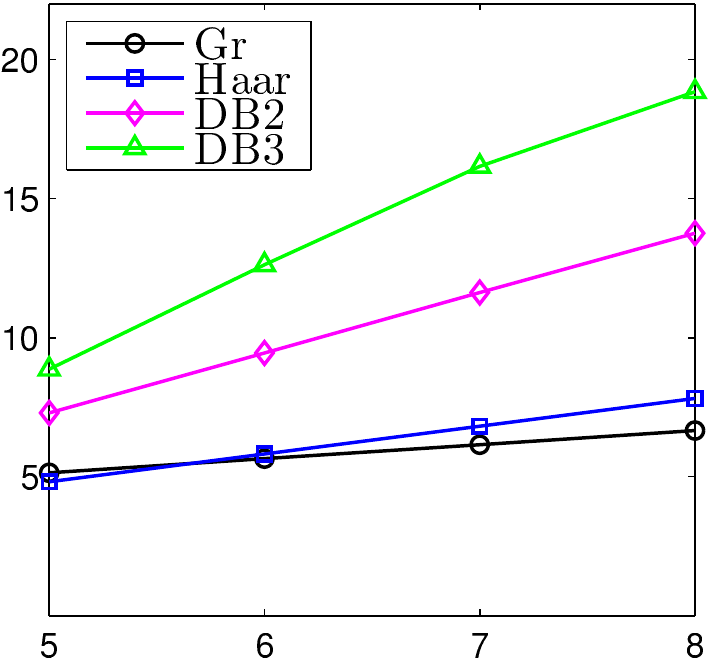} \\
\multicolumn{1}{r}{ $\log M$ \hspace{3.85cm}  $\log K$ \hspace{3.85cm} $\log K$} \vspace{-0.5cm}
\end{tabular}
\end{center}
\caption{\small{
A nonperiodic continuous function $f(x)=x\cos(3\pi x)\chi_{[0,1]}(x)$ is reconstructed from pointwise samples of its Fourier transform taken on an equispaced grid with $\epsilon=1$ (left), on a jittered scheme with jitter $0.1$ (middle) and on a log scheme (right), where the last two satisfy $\delta<0.97$. Reconstruction is performed via GS using different types of boundary-corrected Daubechies wavelets: Haar, DB2 and DB3, and also via truncated Fourier series (TFS) in the uniform case, and gridding in the nonuniform cases, with density compensation factors defined as in \R{eq:weights}.
}}
\label{fig:IncreaseWav1D}
\end{figure}

These superior convergence rates given by Corollary \ref{cor:rates} are depicted in Figure \ref{fig:IncreaseWav1D}. Namely, using  an example of a continuous, nonperiodic function, we compare the convergence rates of generalized sampling with boundary-corrected Daubechies wavelets to the suboptimal convergence rates of the simple direct approaches based on the discretization of the Fourier integral.

\section{Computations in one dimension}\label{s:algorithm_1d}

Let $a$ be the number of vanishing moments of Daubechies wavelets, $R>\log_2(2a)$ the finest wavelet scale and $N=2^R$  the dimension of the reconstruction space $\cW_N$ defined by  \R{eq:rec_space_wavelets}. Let $M\geq N$, and  let $\Omega_M=\left\{\omega_m:m=1,\ldots,M\right\}$ be the given set of sampling points with associated Voronoi weights $\mu_m$ defined as in \R{eq:weights}, and let $\cE_M$ be the sampling space defined by \R{eq:wframe_sampling_space}.
For these choices of the reconstruction and the sampling spaces, generalized sampling reduces to the weighted least-squares system
\be{\label{lsqr_final}
G^{[N,M]}(\alpha) = (\sqrt{\mu_m}\ip{f}{e_{\omega_m}})_{j=1}^M,
}
for a given vector of Fourier samples $(\ip{f}{e_{\omega_m}})_{j=1}^M$ and solving for an unknown vector of boundary-corrected wavelet coefficients $\alpha=(\alpha_n)_{n=1}^N$ of a function $f$.

\rem{As explained previously, efficient implementation of generalized sampling leans on the efficient implementation of the forward and adjoint operations, $G^{[N,M]}$ and $(G^{[N,M]})^*$, which we describe in detail below. They can be implemented as a  function handle in \textsc{Matlab}, for example, and then used in \textsc{Matlab}'s function \texttt{lsqr} for iterative solving of the least squares system \R{lsqr_final}.}

From \R{eq:gs_op}, it follows that, for spaces $\cW_N$ and $\cE_M$, given  $\alpha \in \bbC^{N}$ and $ \zeta \in \bbC^M$, the forward operation can be written as 
\be{\label{forw_1d}
\beta := G^{[N,M]}(\alpha) =  \left(\langle \sum_{n=0}^{N-1} \alpha_n \varphi_n ,\sqrt{\mu_{\omega_m}} e_{\omega_m} \rangle \right)_{m=1}^M, 
}
and the adjoint operation can be written as
\be{\label{adj_1d}
\gamma := (G^{[N,M]})^*(\zeta) = \left( \langle \sum_{m=1}^M \sqrt{\mu_{\omega_m}} \zeta_m e_{\omega_m} , \varphi_n \rangle \right)_{n=0}^{N-1},
} 
We now describe how these operations can be computed efficiently by using the following operators:
\begin{enumerate}
 \item[$i)$] For the set of frequencies $\Omega_M$ and the corresponding set of Voronoi weights $\mu_{m}>0$, the  diagonal weighting operator $V=V_{\Omega_M}:\bbC^{M}\rightarrow\bbC^{M}$ is given by
\be{\label{diag_w}
V(\gamma) =\left( \sqrt{\mu_{m}} \gamma_m \right)_{m=1}^M, \qquad \gamma\in\bbC^{M}.
}
 \item[$ii)$] For the set of frequencies $\Omega_M$, the operator  $F=F_{\Omega_M} : \bbC^{N} \rightarrow \bbC^M$ is given by 
\be{\label{fourier_forw}
F (\gamma)= \left ( \frac{1}{\sqrt{N}} \sum_{k=a}^{N-a-1} \gamma_k  e_{\omega_m} \left(-\frac{k}{N}\right) \right)_{m=1}^{M}, \qquad \gamma\in\bbC^{N}.
}
 \item[$iii)$] For the set of frequencies  $\Omega_M$ and the scaling function $\phi$, the operator  $D=D_{\Omega_M, \phi}: \bbC^M \rightarrow \bbC^M$ is given by 
\be{\label{diag_forw}
D (\zeta) = \left(\hat \phi\left(\frac{\omega_m}{N}\right) \zeta_m \right)_{m=1}^M, \qquad \zeta\in\bbC^{M}.
}
\end{enumerate}
For the weighting operator we have $V^*=V$. The adjoint operator of $F$ is $F^*:\bbC^M \rightarrow \bbC^N$ given by
\be{\label{fourier_adj}
\left( F^* (\zeta) \right)_k=  \begin{cases}
\frac{1}{\sqrt{N}} \sum_{m=1}^{M} \zeta_m  e_{\omega_m} \left(\frac{k}{N}\right)  &k=a,\ldots, N-a-1\\
0 &\text{otherwise}
\end{cases}, \qquad \zeta \in \bbC^M 
}
and the adjoint operator of $D$ is $D^*:\bbC^M \rightarrow \bbC^M$ given by
\be{\label{diag_adj}
D^* (\zeta) =  \left( \overline{ \hat \phi\left(\frac{\omega_m}{N}\right) } \zeta_m \right)_{m=1}^M,\qquad \zeta \in \bbC^M.
}

Now we can analyse the operations \R{forw_1d} and \R{adj_1d}.
We first consider the forward operation. Given $\alpha\in\bbC^N$, by the definition of discrete IWT,  the equation $\beta = G^{[N,M]}(\alpha)$ is equivalent to
$$
\beta_m = \sqrt{\mu_m} \sum_{k=0}^{N-1} \tilde \alpha_k \langle  \phi^{[0,1]}_{R,k} , e_{\omega_m} \rangle,   \quad m=1,\ldots,M,
$$
where $\tilde \alpha =W^{-1}(\alpha) \in \bbC^{N}$.
Since the Fourier transform of the internal scaling function $\langle \phi^{[0,1]}_{R,k}, e_{\omega}  \rangle$  can be written as
$$
\hat{\phi}^{[0,1]}_{R,k}(\omega) = \frac{1}{\sqrt{N}}\hat{\phi}\left(\frac{\omega}{N}\right) e_\omega \left(-\frac{k}{N}\right), \quad k=a,\ldots,N-a-1,
$$
by using definitions of operators $F$ and $D$, we get
$$
\tilde\beta_m=\frac{1}{\sqrt{\mu_m}}\beta_m = \frac{1}{\sqrt{N}}\sum_{k=0}^{a-1} \tilde \alpha_k \hat{\phi}^0_k\left(\frac{\omega_m}{N}\right) + \left(D \left( F (\tilde\alpha) \right) \right)_m \!\! + \frac{1}{\sqrt{N}}\!\! \sum_{k=N-a}^{N-1}\!\!\! \tilde \alpha_k \hat{\phi}^1_{N-k-1}\left(\frac{\omega_m}{N}\right), \ m=1,\ldots,M.
$$
Once $\tilde \beta$ has been computed, it is left to apply the weighting operator and get $\beta=V(\tilde\beta)$.

For the adjoint operation, to compute $\gamma=\left( G^{[N,M]}\right)^* (\zeta)$  given $\zeta\in\bbC^{M}$,  we first apply the  weighting operator and set $\tilde \zeta = V(\zeta)$. Then, similarly to the forward operation case, one can check that $\tilde \gamma=W^{-1}\gamma$ and $\zeta$ are related by the following equations.
$$
\tilde \gamma_k = \frac{1}{\sqrt{N}}  \sum_{m=1}^{M} \tilde\zeta_m \overline{\hat{\phi}^0_k\left(\frac{\omega_m}{N}\right)} , \ k=0,\ldots,a-1, \quad \tilde \gamma_k = \frac{1}{\sqrt{N}} \sum_{m=1}^{M} \tilde\zeta_m \overline{\hat{\phi}^1_{N-k-1}\left(\frac{\omega_m}{N}\right)}, \ k=N-a, \ldots, N-1, 
$$
and 
$$
\tilde \gamma_k = \frac{1}{\sqrt{N}} \sum_{m=1}^M \overline{\hat{\phi}\left(\frac{\omega_m}{N}\right)} \tilde\zeta_m  e_{\omega_m}\left(\frac{k}{N}\right), \quad k=a,\ldots,N-a-1.
$$
Note that, by using adjoint operators $D^*$ and $F^*$, this last part can be written as
$$
\tilde \gamma_k =  \left(F^* \left ( D^*(\tilde\zeta) \right)\right)_k, \quad k=a,\ldots,N-a-1. 
$$

These computational steps, which we summarize below, lead to the efficient algorithm for forward and adjoint operations, and therefore to to the efficient algorithm for solving the weighted least-squares system \R{lsqr_final}.

\subsection*{The one-dimensional algorithm}\label{ss:1algorithm}

Precompute the Voronoi weights $(\mu_m)_{m=1}^M$ and pointwise measurements of the Fourier transforms of the three scaling functions: 
$$
\left(\hat \phi\left(\frac{\omega_m}{N}\right)\right)_{m=1}^M, \quad \left(\hat{\phi}_k^0\left(\frac{\omega_m}{N}\right)\right)_{m=1}^M, \quad \left(\hat{\phi}_k^1\left(\frac{\omega_m}{N}\right)\right)_{m=1}^M,\quad k=1,\ldots, a.
$$

\subsubsection*{The forward operation}
 
Given $\alpha \in \bbC^{N}$, $\beta = G^{[N,M]}(\alpha)$ can be obtained by applying the following steps.
\begin{enumerate}
\item Compute the scaling coefficients $\tilde \alpha = W^{-1} (\alpha)$, where $W^{-1}$ is the one-dimensional discrete boundary-corrected IWT.
\item Compute contributions from the boundary scaling functions:
$$
\tilde\beta^0 = \left( \frac{1}{\sqrt{N}} \sum_{k=0}^{a-1} \tilde \alpha_k \hat{\phi}^0_k\left(\frac{\omega_m}{N}\right) \right)_{m=1}^M, \quad  \tilde\beta^1 = \left( \frac{1}{\sqrt{N}} \sum_{k=N-a}^{N-1} \tilde \alpha_k \hat{\phi}^1_{N-k-1}\left(\frac{\omega_m}{N}\right)\right)_{m=1}^M.
$$
\item  Compute contribution from the internal scaling functions:
\begin{enumerate}[label*=\arabic*.]
 \item Apply $F$ to $\tilde \alpha$ to get
$\hat \alpha = F\left(  \tilde \alpha \right)$, where $F$ is defined by \R{fourier_forw}.
 \item Apply $D$ to $\hat \alpha$ to get
$\tilde\beta^{\text{int}} = D(\hat \alpha)$, where $D$ is defined by \R{diag_forw}.
\end{enumerate}
\item Compute $\tilde\beta =  \tilde\beta^0 + \tilde\beta^1 + \tilde\beta^{\text{int}}$.
\item Apply $V$ to compute $\beta=V(\tilde\beta)$, where $V$ is defined by \R{diag_w}.
\end{enumerate}

\subsubsection*{The adjoint operation}

Given $\zeta \in \bbC^{M}$, $\gamma = \left(G^{[N,M]}\right)^*(\zeta)$ can be computed as follows.

\begin{enumerate}
 \item  Apply the weighting operator $V$ and set $\tilde \zeta = V(\zeta)$.
\item Compute the coefficients of the boundary scaling functions:
$$
\tilde \gamma_k = \frac{1}{\sqrt{N}}  \sum_{m=1}^{M} \tilde\zeta_m \overline{\hat{\phi}^0_k\left(\frac{\omega_m}{N}\right)} ,  \quad \tilde \gamma_{N-k-1} = \frac{1}{\sqrt{N}} \sum_{m=1}^{M} \tilde\zeta_m \overline{\hat\phi^1_{k}\left(\frac{\omega_m}{N}\right)},  
$$
for $k=0,\ldots,a-1$.
\item Compute the coefficients of the internal scaling functions:
\begin{enumerate}[label*=\arabic*.]
 \item Compute $\tilde\zeta_{\phi} = D^*(\tilde\zeta)$, where $D^*$ is defined by \R{diag_adj}.
 \item Compute $\tilde \gamma_k = \left(F^*(\tilde\zeta_{\phi})\right)_k$, $k=a-1,\ldots,N-a-1$, where $F^*$ is defined by \R{fourier_adj}.
\end{enumerate}
\item Compute $\gamma = W(\tilde \gamma)$, where $W$ is a discrete one-dimensional boundary-corrected FWT.\vspace{0.3cm}
\end{enumerate}

\rem{Note that when the sampling points are uniform with spacing $\epsilon$, the weights are simply $\mu_m=\epsilon$ for all $m$, and the precomputation of Voronoi weights is not required. }

\rem{The above algorithm requires the precomputation of pointwise evaluations of the Fourier transform of the internal and boundary scaling functions. Note that for the internal scaling function $\phi$, we may use the approximation $$
  \prod_{j=1}^N m_0(2^{-j} \xi) \to \hat \phi(\xi), \qquad N\to \infty
$$
where $m_0$ is a trigonometric polynomial \cite{daubechies1992ten}. A similar approximation may be used in the case of the boundary scaling functions. This is outlined in detail in the appendix.
}

\rem{
On evaluating the reconstructed signal: Recall that in solving 
\bes{
G^{[N,M]}(\alpha) = (\ip{f}{\psi_j})_{j=1}^M.
}
for an appropriate $M=\ord{N}$, we obtain $\alpha$ which is an approximation of the first $N$ wavelet coefficients of $f$ and the reconstructed signal is $\tilde f = \sum_{j=1}^N \alpha_j \phi_j$. To evaluate $\tilde f$ on the grid points $(k 2^{-J})_{k=1}^{2^J}$ for $J\in\bbN$, it suffices to evaluate each $\phi$ on these grid points and we may do so by either implementing the cascade algorithm \cite{daubechies1992ten} or the dyadic dilation algorithm \cite{louis1997wavelets}.
}

\paragraph{Computational cost of the one-dimensional algorithm.}
Let us analyse the computational cost of the forward operation. The adjoint operation can be analysed similarly leading to the same computational cost. The computational cost of applying the discrete boundary-corrected IWT is $\ord{N}$. The cost of step 2, involving boundary scaling functions, is $\ord{aM}$. For step 3.1, the key point is to observe that $F$ is simply a restricted and shifted version of the discrete (nonuniform) Fourier transform,  and thus its fast implementation (FFT/NUFFT) can be used when computing $F(\tilde \alpha)$. Hence, the the cost of step 3.1 is $\ord{M\log(N/\epsilon)}$ in the uniform case and $\ord{L\log(N)+JM}$ in the nonuniform case, where $L$ is the length of underlying interpolating FFT for NUFFT, and $J$ is the number of interpolating coefficients (typically $J=7$) \cite{FesslerNUFFT}. Finally, the cost of the diagonal operations in both steps 3.2 and 5 is $\ord{M}$. Therefore, given that $\epsilon \sim 1$, $J\sim a$ and $L\sim M$, the total cost is essentially $\ord{aM+M\log(N)}$.

\section{Computations in two dimensions}\label{s:algorithm_2d}

Again, let $a$ be the number of vanishing moments, $R>\log_2(2a)$ the finest wavelet scale and  $N=2^R$. Let $\cW_{N^2}$  be the reconstruction space defined by \R{eq:rec_space_wavelets_2d}. Let $M\geq N$, and let $$\Omega_M=\left\{\omega_m\ :\ \omega_m=(\omega_m^1,\omega_m^2),\ m=1,\ldots,M\right\}$$ be the set of sampling points in $\bbR^2$, which we write as $\Omega_M=(\Omega_M^1,\Omega_M^2)$ correspondingly. Let $\cE_M$ be the associated sampling space defined by \R{eq:wframe_sampling_space}, with the Voronoi weights $\mu_m$ defined as in \R{eq:weights}.
In this case, the least-squares system \R{eq:least_sq}  becomes 
$$
\left(  \mu_m \ip{\sum_{n_1,n_2=1}^N \alpha_{n_1,n_2} \varphi_{n_1,n_2}}{e_{\omega_m}}  \right)_{m=1}^M = \left( \mu_m \ip{f}{e_{\omega_m}} \right)_{m=1}^{M}.
$$
If we apply the two-dimensional boundary-corrected IWT to the matrix of wavelet coefficients $\alpha\in\bbC^{N\times N}$, so that $\tilde{\alpha}=\textit{\textbf{W}}^{-1}(\alpha)\in\bbC^{N\times N}$, we get
$$
\left( \mu_m  \ip{\sum_{k_1,k_2=1}^N \tilde\alpha_{k_1,k_2} \Phi_{R,(k_1,k_2)}}{e_{\omega_m}}  \right)_{m=1}^M = \left(\mu_m\ip{f}{e_{\omega_m}} \right)_{m=1}^{M}.
$$
Since $\Phi_{R,(k_1,k_2)}(x_1,x_2)=\phi_{R,k_1}^{[0,1]}(x_1) \phi_{R,k_2}^{[0,1]}(x_2)$, we can write the following algorithm.

\subsection*{The two-dimensional algorithm} 

Precompute the vectors $(\mu_m)_{m=1}^M$ and
$$
\left(\hat \phi\left(\frac{\omega^{i}_m}{N}\right)\right)_{m=1}^M, \quad \left(\hat{\phi}_k^0\left(\frac{\omega^{i}_m}{N}\right)\right)_{m=1}^M, \quad \left(\hat{\phi}_k^1\left(\frac{\omega^{i}_m}{N}\right)\right)_{m=1}^M,\quad k=1,\ldots, a, \quad  i=1,2.
$$

\subsubsection*{The forward operation}
 
Given $\alpha \in \bbC^{N\times N}$, $\beta = G^{[N^2,M]}(\alpha)\in\bbC^M$ can be obtained by applying the following steps.
\begin{enumerate}
\item Compute the scaling coefficients $\tilde \alpha = \textit{\textbf{W}}^{-1} (\alpha)$, where $\textit{\textbf{W}}^{-1}$ is the discrete two-dimensional boundary corrected IWT.
\item Compute contributions from the corners (the boundary scaling functions in the both axes):
\eas{
\beta^{00} &= \left( \frac{1}{N}  \sum_{k_1=0}^{a-1}  \sum_{k_2=0}^{a-1} \tilde\alpha_{k_1,k_2} \hat\phi_{k_1}^0 \left(\frac{\omega_m^1}{N}\right) \hat\phi_{k_2}^0\left(\frac{\omega_m^2}{N}\right) \right)_{m=1}^M \\
\beta^{01} &= \left( \frac{1}{N}  \sum_{k_1=0}^{a-1} \sum_{k_2=N-a}^{N-1} \tilde\alpha_{k_1,k_2}  \hat\phi_{k_1}^0 \left(\frac{\omega_m^1}{N}\right) \hat\phi_{N-k_2-1}^1\left(\frac{\omega_m^2}{N}\right)  \right)_{m=1}^M \\
\beta^{10} &= \left( \frac{1}{N}  \sum_{k_1=N-a}^{N-1}  \sum_{k_2=0}^{a-1} \tilde\alpha_{k_1,k_2} \hat\phi_{N-k_1-1}^1 \left(\frac{\omega_m^1}{N}\right) \hat\phi_{k_2}^0\left(\frac{\omega_m^2}{N}\right)  \right)_{m=1}^M \\
\beta^{11} &= \left( \frac{1}{N} \sum_{k_1=N-a}^{N-1}  \sum_{k_2=N-a}^{N-1} \tilde\alpha_{k_1,k_2} \hat\phi_{N-k_1-1}^1 \left(\frac{\omega_m^1}{N}\right)  \hat\phi_{N-k_2-1}^1\left(\frac{\omega_m^2}{N}\right)  \right)_{m=1}^M 
}
\item Compute contributions from the edges (the boundary scaling functions in only one of the axes):
\eas{
\tilde \beta^{0\text{int}} &= \frac{1}{\sqrt{N}}  \sum_{k_1=0}^{a-1}  D_{\Omega_M^1 \phi_{k_1}^0}  D_{\Omega_M^2 \phi} F_{\Omega_M^2}\left(\tilde{\alpha}_{k_1,\cdot}\right)  \\
\tilde \beta^{1\text{int}} &= \frac{1}{\sqrt{N}} \sum_{k_1=N-a}^{N-1}  D_{\Omega_M^1 \phi_{N-k_1-1}^1}  D_{\Omega_M^2 \phi} F_{\Omega_M^2}\left(\tilde{\alpha}_{k_1,\cdot}\right) \\
\tilde \beta^{\text{int}0} &= \frac{1}{\sqrt{N}} \sum_{k_2=0}^{a-1}  D_{\Omega_M^1 \phi }  D_{\Omega_M^2 \phi_{k_2}^0} F_{\Omega_M^1}\left(\tilde{\alpha}_{\cdot,k_2}\right)  \\
\tilde \beta^{\text{int}1} &= \frac{1}{\sqrt{N}} \sum_{k_2=0}^{a-1}  D_{\Omega_M^1 \phi }  D_{\Omega_M^2 \phi_{N-k_2-1}^1} F_{\Omega_M^1}\left(\tilde{\alpha}_{\cdot,k_2}\right)  
}
where $F$ and $D$ are defined by \R{fourier_forw} and \R{diag_forw}, respectively.
\item  Compute the contribution from the internal scaling functions:
\begin{enumerate}[label*=\arabic*.]
 \item $\hat \alpha = \textbf{F}_{\Omega_M} \left( \tilde{\alpha} \right)$, where $\textbf{F}_{\Omega_M}:\bbC^{N\times N} \rightarrow \bbC^M$ is such that for each $\gamma\in\bbC^{N\times N}$
$$
\textbf{F}_{\Omega_M} (\gamma) = \left( \frac{1}{N} \sum_{k_1,k_2=a}^{N-a-1} \gamma_{k_1,k_2} e_{\omega_m} \left(-\frac{(k_1,k_2)}{N}\right) \right)_{m=1}^{M}.
$$
 \item $\tilde \beta^{\text{int}\text{int}} = D_{\Omega_M^1, \phi} D_{\Omega_M^2, \phi}(\hat \alpha)$.
\end{enumerate}
\item Compute $\tilde \beta =  \tilde \beta^{00}+\tilde \beta^{01}+\tilde \beta^{10}+\tilde \beta^{11} + \tilde \beta^{0\text{int}} + \tilde \beta^{1\text{int}} + \tilde \beta^{\text{int}0} + \tilde \beta^{\text{int}1}  + \tilde \beta^{\text{int}\text{int}}$.
\item Apply $V$ to get $\beta=V(\tilde\beta)$, where $V$ is defined by \R{diag_w}.
\end{enumerate} 

\subsubsection*{The adjoint operation}

Given $\zeta \in \bbC^{M}$, $\gamma = \left(G^{[N^2,M]}\right)^*(\zeta)\in\bbC^{N,N}$ can be computed as follows.

\begin{enumerate}
 \item Apply the weighting operator $V$ and set $\tilde\zeta=V(\zeta)$.
\item Compute the scaling coefficients at the corners
\eas{
\tilde\gamma_{k_1,k_2} &= \frac{1}{N} \sum_{m=1}^{M} \zeta_m \overline{\hat\phi^0_{k_1} \left(\frac{\omega_m^1}{N}\right)} \overline{\hat\phi^0_{k_2} \left(\frac{\omega_m^2}{N}\right)}, \\
\tilde\gamma_{k_1,N-a+k_2} &= \frac{1}{N} \sum_{m=1}^{M} \zeta_m \overline{\hat\phi^0_{k_1} \left(\frac{\omega_m^1}{N}\right)} \overline{\hat\phi^1_{a-k_2-1} \left(\frac{\omega_m^2}{N}\right)}\\
\tilde\gamma_{N-a+k_1,k_2} &= \frac{1}{N} \sum_{m=1}^{M} \zeta_m \overline{\hat\phi^1_{a-k_1-1} \left(\frac{\omega_m^1}{N}\right)} \overline{\hat\phi^0_{k_2} \left(\frac{\omega_m^2}{N}\right)}, \\ 
\tilde\gamma_{N-a+k_1,N-a+k_2} &= \frac{1}{N} \sum_{m=1}^{M} \zeta_m \overline{\hat\phi^1_{a-k_1-1} \left(\frac{\omega_m^1}{N}\right)} \overline{\hat\phi^1_{a-k_2-1} \left(\frac{\omega_m^2}{N}\right)}.
}
for $k_1,k_2=0,\ldots,a-1$.
\item Compute the scaling coefficients at the edges
\eas{
\tilde \gamma_{k_1,k_2} &= \frac{1}{\sqrt{N}} \left( \left( F_{\Omega_M^1} \right)^* \left( D_{\Omega_M^1 \phi}\right)^* \left( D_{\Omega_M^2 \phi^0_{k_2}} \right)^* (\tilde\zeta) \right)_{k_1},\\
\tilde \gamma_{k_1,k_2} &= \frac{1}{\sqrt{N}} \left( \left( F_{\Omega_M^1} \right)^* \left( D_{\Omega_M^1 \phi}\right)^* \left( D_{\Omega_M^2 \phi^1_{a-k_2-1}} \right)^* (\tilde\zeta) \right)_{k_1},
}
for $k_1=a,\ldots,N-a-1, \ k_2=0,\ldots,a-1$ and
\eas{
\tilde \gamma_{k_1,k_2} &= \frac{1}{\sqrt{N}} \left( \left( F_{\Omega_M^2} \right)^* \left( D_{\Omega_M^1 \phi^0_{k_1}}\right)^* \left( D_{\Omega_M^2 \phi} \right)^* (\tilde\zeta) \right)_{k_2 },\\
\tilde \gamma_{k_1,k_2} &= \frac{1}{\sqrt{N}} \left( \left( F_{\Omega_M^2} \right)^* \left( D_{\Omega_M^1 \phi^1_{a-k_1-1}}\right)^* \left( D_{\Omega_M^2 \phi} \right)^* (\tilde\zeta) \right)_{k_2 }, 
}
for $k_1=0,\ldots,a-1, \ k_2=a,\ldots,N-a-1$, where $F^*$ and $D^*$ are defined by \R{fourier_adj} and \R{diag_adj}, respectively.
\item Compute the scaling coefficients of the internal wavelets
\begin{enumerate}[label*=\arabic*.]
 \item $\tilde\zeta_{\phi,\phi} = \left(D_{\Omega_M^1, \phi}\right)^* \left(D_{\Omega_M^2, \phi}\right)^*(\tilde\zeta)$.
 \item $\tilde \gamma_{k_1,k_2} = \left(\textbf{F}^*\left(\tilde\zeta_{\phi,\phi}\right)\right)_{k_1,k_2}$, \ $k_1,k_2=a-1,\ldots,N-a-1$.
\end{enumerate}
\item Compute $\gamma = \textit{\textbf{W}}(\tilde \gamma)$, where $\textit{\textbf{W}}$ is the discrete two-dimensional boundary corrected FWT.
\end{enumerate}

\paragraph{Computational cost of the two-dimensional algorithm.}
For the forward operation, the computational cost of  step 1 is $\ord{N^2}$ and that of step 2 is $\ord{a^2M}$. Step 3 is $\ord{a(M+M\log(N/\epsilon))}$ in the uniform case and $\ord{a(M+L\log(N)+JM)}$ computations in the nonuniform case. The cost of step 4.1 is basically the cost of the two-dimensional NUFFT or FFT, i.e, $\ord{L^2\log(N^2)+J^2M}$ or $\ord{M\log(N^2/\epsilon^2)}$. The cost of step 4.2  as well as step 6 is $\ord{M}$. Hence, if we assume $\epsilon \sim 1$, $J\sim a$ and $L^2\sim M$, the total cost is $\ord{a^2 M+M\log N^2}$. The same cost holds for the adjoint operation.

\section{Simplifications in the uniform  two-dimensional case}\label{s:algorithm_2d_simp}

Let us assume the setting from the previous section, but now the sampling points are given uniformly i.e, the sampling space $\cE_{M^2}$ is given by \R{sampling_space_uniform2d} where, for convenience and without loss of generality, we assume $\epsilon_1=\epsilon_2=1$. 
In this case, we have
$$
G^{[N^2,M^2]}: \bbC^{N\times N}\to \bbC^{M\times M},\quad \alpha \mapsto \left( \ip{\sum_{n_1,n_2=1}^N \alpha_{n_1,n_2} \varphi_{n_1,n_2}}{e_{m_1} e_{m_2}}\right)_{-\lceil M/2\rceil \leq m_1, m_2 \leq \lceil M/2\rceil -1}.
$$

Let us analyze the forward operation, i.e. the computation of $\beta = G^{[N^2,M^2]}(\alpha)$, for a given $\alpha$. As usual, we can apply the discrete wavelet transform to get
$$
\beta_{m_1, m_2} = \ip{ \sum_{k_1,k_2=0}^{N-1} \tilde \alpha_{k_1,k_2} \Phi_{R,(k_1,k_2)}}{e_{m_1} e_{m_2}}, \qquad -\lceil  M/2\rceil \leq \abs{m_1},\abs{m_2} \leq M
$$
where $\tilde \alpha =\textit{\textbf{W}}^{-1}\alpha$. Now denote steps 2-3 of the forward operation in the one-dimensional algorithm applied to the scaling coefficients by $G_{\phi}:\bbC^{N}\rightarrow\bbC^M$. Given the special structure of both the reconstruction and the sampling vectors, i.e, by using the definition of  $\Phi_{R,(k_1,k_2)}$ and $e_{m_1} e_{m_2}$, we obtain
\spls{
\beta_{m_1, m_2} = \ip{ \sum_{k_1=0}^{N-1} \ip{  \sum_{k_2=0}^{N-1} \tilde \alpha_{k_1,k_2} \phi_{R,k_2}^{[0,1]}}{e_{ m_2}}\phi_{R,k_1}^{[0,1]}}{e_{m_1}}
= \ip{\sum_{k_1=0}^{N-1} \gamma_{k_1, m_2}\phi_{R,k_1}^{[0,1]}}{e_{m_1}} = \eta_{m_1,m_2}
}
where for each $k_1=0,\ldots,N-1$ 
$$
\left(\gamma_{k_1, m_2}\right)_{m_2=-\lceil  M/2\rceil}^{\lceil  M/2\rceil-1} = \left( \ip{  \sum_{k_2=0}^{N-1} \tilde \alpha_{k_1,k_2} \phi_{R,k_2}^{[0,1]}}{e_{{m_2}}} \right)_{m_2=-\lceil  M/2\rceil}^{\lceil  M/2\rceil-1} = G_{\phi}\left(\left(\tilde \alpha_{k_1,k_2}\right)_{k_2=0}^{N-1}\right)
$$
and for each $m_2=-\lceil  M/2\rceil, \ldots, \lceil  M/2\rceil-1$
$$
\left(\eta_{m_1,m_2} \right)_{m_1=-\lceil  M/2\rceil}^{\lceil  M/2\rceil-1} = G_{\phi}\left(\left(\gamma_{k_1,m_2}\right)_{k_1=0}^{N-1}\right).
$$
Hence we can write the following algorithm.

\subsection*{The two-dimensional algorithm in the uniform case}\label{s:algorithm_2d_uniform}

Precompute the vectors 
$$
\left(\hat \phi\left(\frac{  m}{N}\right)\right)_{m=-\lceil  M/2\rceil}^{\lceil  M/2\rceil-1}, \quad \left(\hat{\phi}_k^0\left(\frac{  m}{N}\right)\right)_{m=-\lceil  M/2\rceil}^{\lceil  M/2\rceil-1}, \quad \left(\hat{\phi}_k^1\left(\frac{  m}{N}\right)\right)_{m=-\lceil  M/2\rceil}^{\lceil  M/2\rceil-1},\quad k=1,\ldots, a.
$$

\subsubsection*{The forward operation}

Given $\alpha\in\bbC^{N\times N}$, $\beta = G^{[N^2,M^2]}(\alpha)$ can be obtained as follows.
\begin{enumerate}
\item Compute the scaling coefficients $\tilde \alpha = \textit{\textbf{W}}^{-1} (\alpha)$. 
\item For each $k_1=0,\ldots, N-1$, apply $G_{\phi}$ to the $n_1^{th}$ row of $\tilde \alpha$ to obtain 
$$
\tilde \beta_{k_1} = G_{\phi}\left( \left(\tilde \alpha_{k_1,k_2}\right)_{k_2=0}^{N-1}\right) \in \bbC^M.
$$
Let 
$$
\tilde \beta = \begin{pmatrix}
\tilde \beta_0^T\\
\vdots\\
\tilde \beta_{N-1}^T
\end{pmatrix}\in\bbC^{N\times M}.
$$
\item For each $m_2=-\lceil M/2\rceil,\ldots, \lceil M/2\rceil-1$, apply $G_{\phi}$ to the $m_2^{th}$ column of $\tilde \beta$ to obtain 
$$
\beta_{m_2} =  G_{\phi}\left( \left(\tilde \beta_{k_1, m_2}\right)_{k_1=0}^{N-1}\right) \in \bbC^M.
$$
Let
$$
\beta = \begin{pmatrix}
\beta_{-\lceil M/2\rceil}|\hdots|\beta_{\lceil M/2\rceil-1}
\end{pmatrix}\in\bbC^{M\times M}.
$$
\end{enumerate} 

\subsubsection*{The adjoint operation}
Given $\zeta\in\bbC^{M\times M}$, $\gamma=\left(G^{[N^2,M^2]}\right)^*(\zeta)$ can be computed as follows.
\begin{enumerate}
\item For each $k_2=-\lceil M/2\rceil,\ldots, \lceil M/2\rceil-1$, let the column of $\zeta$ indexed by $k_2$ be denoted by $\zeta^{k_2}$, so $\zeta^{k_2} = \left(\zeta_{k_1,k_2}\right)_{k_1=-\lceil M/2\rceil}^{\lceil M/2\rceil-1}$. Apply $\left(G_\phi\right)^*$ to every column of $\zeta$ to obtain
$$
\gamma = \begin{pmatrix}
         \left(G_\phi\right)^*(\zeta^{-\lceil M/2\rceil}) & \left(G_\phi\right)^*(\zeta^{-\lceil M/2\rceil+1}) & \hdots & \left(G_\phi\right)^*(\zeta^{\lceil M/2\rceil-1}) 
        \end{pmatrix}\in\bbC^{N\times M}.
$$
\item For $n_1=0,\ldots, N-1$, let the row of $\gamma$ indexed by $n_1$ be denoted by $\gamma^{n_1}$. Apply $\left(G_\phi\right)^*$ to each $\gamma^{n_1}$ to obtain
$$
\tilde \gamma = \begin{pmatrix}
\left(\left(G_\phi\right)^*(\gamma^0)\right)^T\\
\vdots\\
\left(\left(G_\phi\right)^*(\gamma^{N-1})\right)^T
\end{pmatrix}\in \bbC^{N\times N}.
$$
\item $\gamma= \textit{\textbf{W}}(\tilde \gamma)$.
\end{enumerate}

\paragraph{Computational cost of the simplified two-dimensional algorithm.}
Again, for the first step of the forward operation we have the computational cost $\ord{N^2}$. Since the cost of each application of operation $G_{\phi}$ is $\ord{aM+M\log(N/\epsilon)}$, and the cost of each transposing of an $M$-dimensional vector is $\ord{M}$, we have that the  cost of step 2 is $\ord{N(aM+M\log(N/\epsilon)+M)}$ and the cost of step 3 is $\ord{M(aM+M\log(N/\epsilon)+M)}$. Therefore, the total  cost amounts to $\ord{aM^2 +M^2\log(N/\epsilon)}$. Note that the cost of the non-simplified two-dimensional algorithm in this case is $\ord{a^2M^2+M^2\log (N/\epsilon)^2}$.

\section{Numerical examples}\label{s:num_ex}

Finally, we present a few numerical examples  using the algorithms from Sections \ref{s:algorithm_1d}-\ref{s:algorithm_2d_uniform}. MATLAB scripts for reproducing all the examples in this section are available at  \url{http://www.damtp.cam.ac.uk/research/afha/code/}.

\begin{figure}[H]
\begin{tabular}{cccc}
                    & Equispaced sampling & Jittered sampling & Log sampling \vspace{0.15cm} \\
{\rotatebox{90}{\hspace{1cm}GS with DB4}} & \includegraphics[scale=0.62]{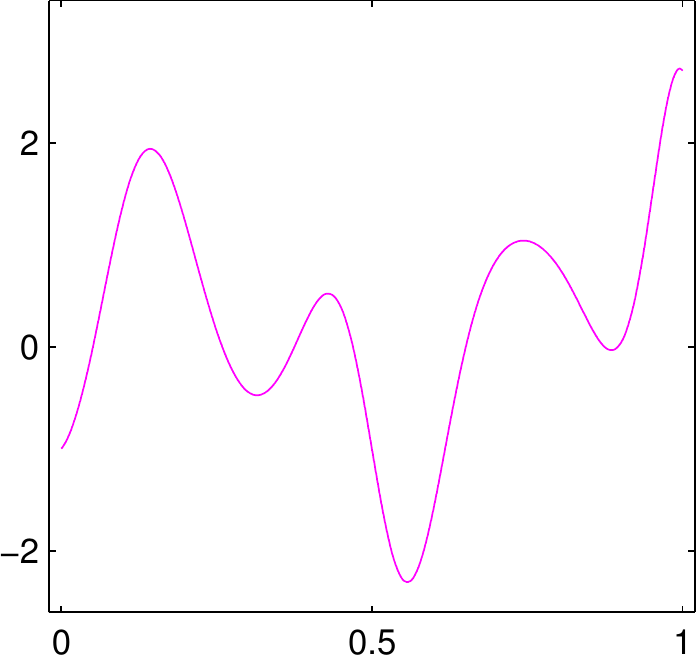} & \includegraphics[scale=0.62]{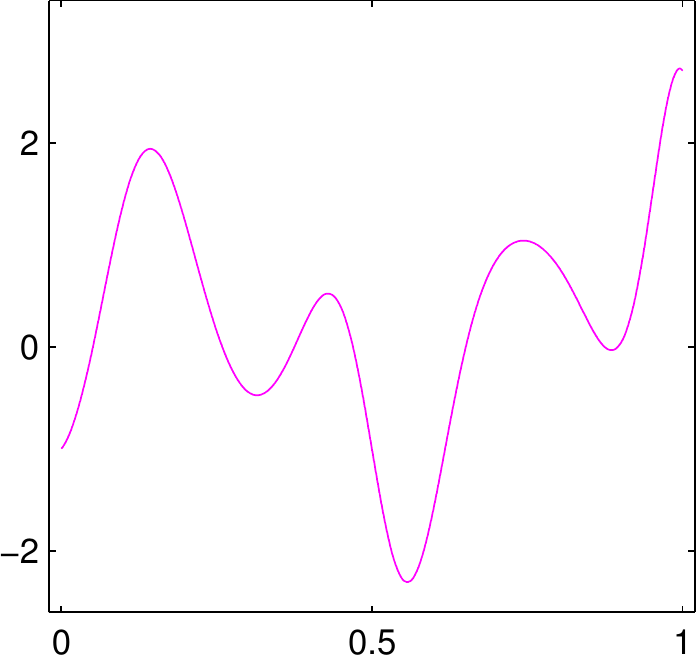} & \includegraphics[scale=0.62]{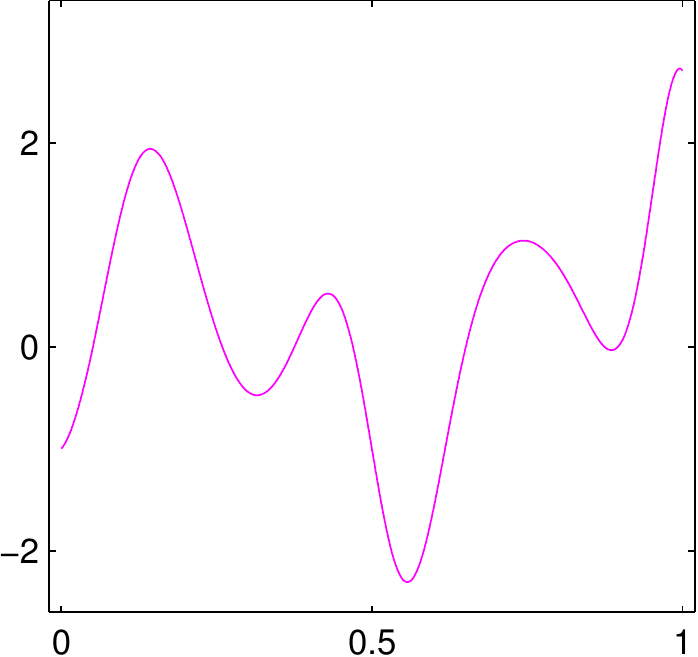} \\
                    & \small{$\qquad5.78\times10^{-4}$} \vspace{0.15cm} & \small{$\qquad5.57\times10^{-4}$}  & \small{$\qquad5.58\times10^{-4}$} \\                 
{\rotatebox{90}{\hspace{1.4cm}TFS/gridding}} & \includegraphics[scale=0.62]{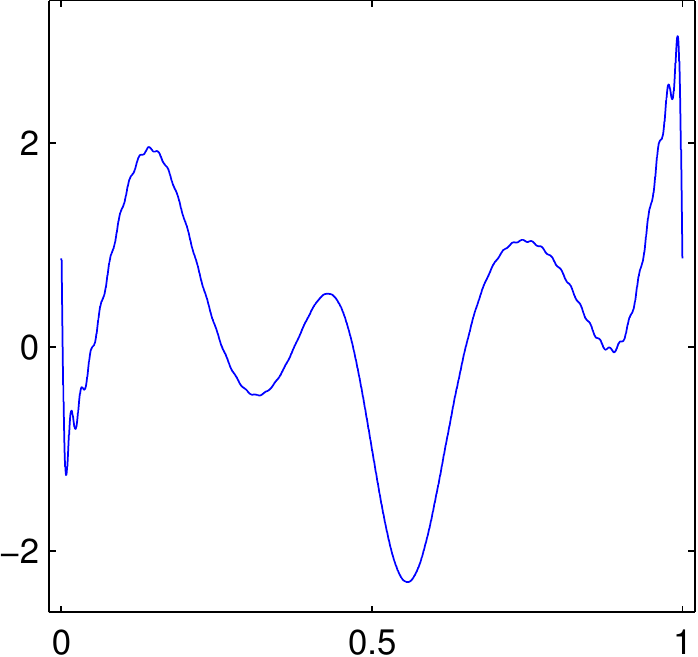} & \includegraphics[scale=0.62]{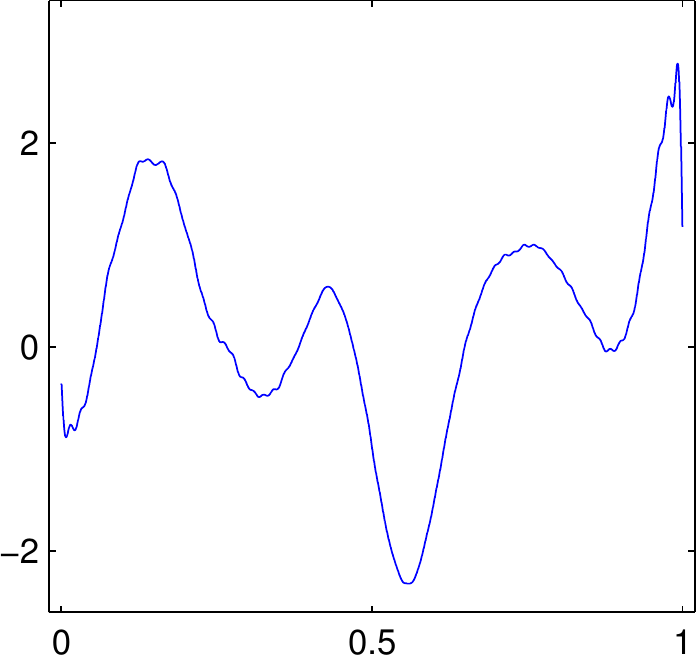} & \includegraphics[scale=0.62]{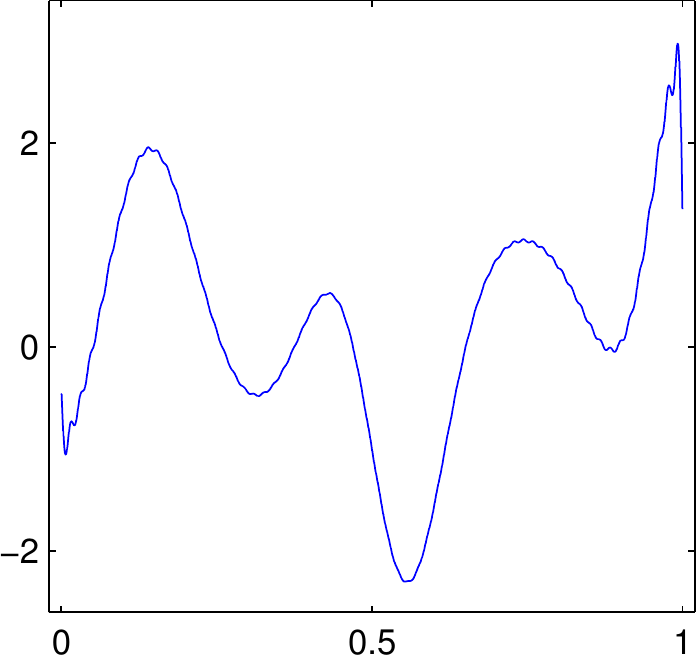} \\
                    & \small{$\qquad1.05\times10^{-1}$} & \small{$\qquad1.27\times10^{-1}$} & \small{$\qquad5.98\times10^{-2}$}                                            
\end{tabular}
\caption{\small{Function \R{ex:fun} is reconstructed from its Fourier samples taken from the interval $[-64,64]$ according to three different sampling schemes: equispaced (with spacing $\epsilon=1$, in total $128$ points), jittered (with spacing $\epsilon=0.77$ and jitter $0.1$, in total $168$ points) and log (with density $\delta=0.97$, in total $653$ points). In the upper panels, reconstruction is done via GS with $64$ DB4 wavelets, while in the lower panels via truncated Fourier series (TFS) in the uniform case, and in the nonuniform cases via gridding with density compensation factors \R{eq:weights}. The $\cL_2$ error is written below each reconstruction.
}}
\label{fig:1Dcompare}
\end{figure}

\paragraph{Example 1.}
In Figure \ref{fig:1Dcompare}, we reconstruct a one-dimensional continuous, but nonperiodic, function 
\be{\label{ex:fun}
f(x)=-\exp(x\cos(4\pi x))\cos(7\pi x)\chi_{[0,1]}(x)+\sin(3\pi x)\chi_{[0,1]}(x)
}
from its Fourier samples taken on three different sampling schemes: equispaced, jittered and log. We use GS with boundary-corrected Daubechies wavelets of order 4 (DB4), and compare its performance to the direct approaches based on a discretization of the Fourier integral (truncated Fourier series and gridding).  While the direct approaches are apparently plagued by Gibbs artifacts and the reconstruction quality evidently depends on the underlying sampling scheme, GS performs equally well for all three sampling schemes producing a numerical error of order $10^{-4}$ by using only $64$ reconstruction functions.

\begin{figure}
\begin{center}
\includegraphics[width=0.3\textwidth]{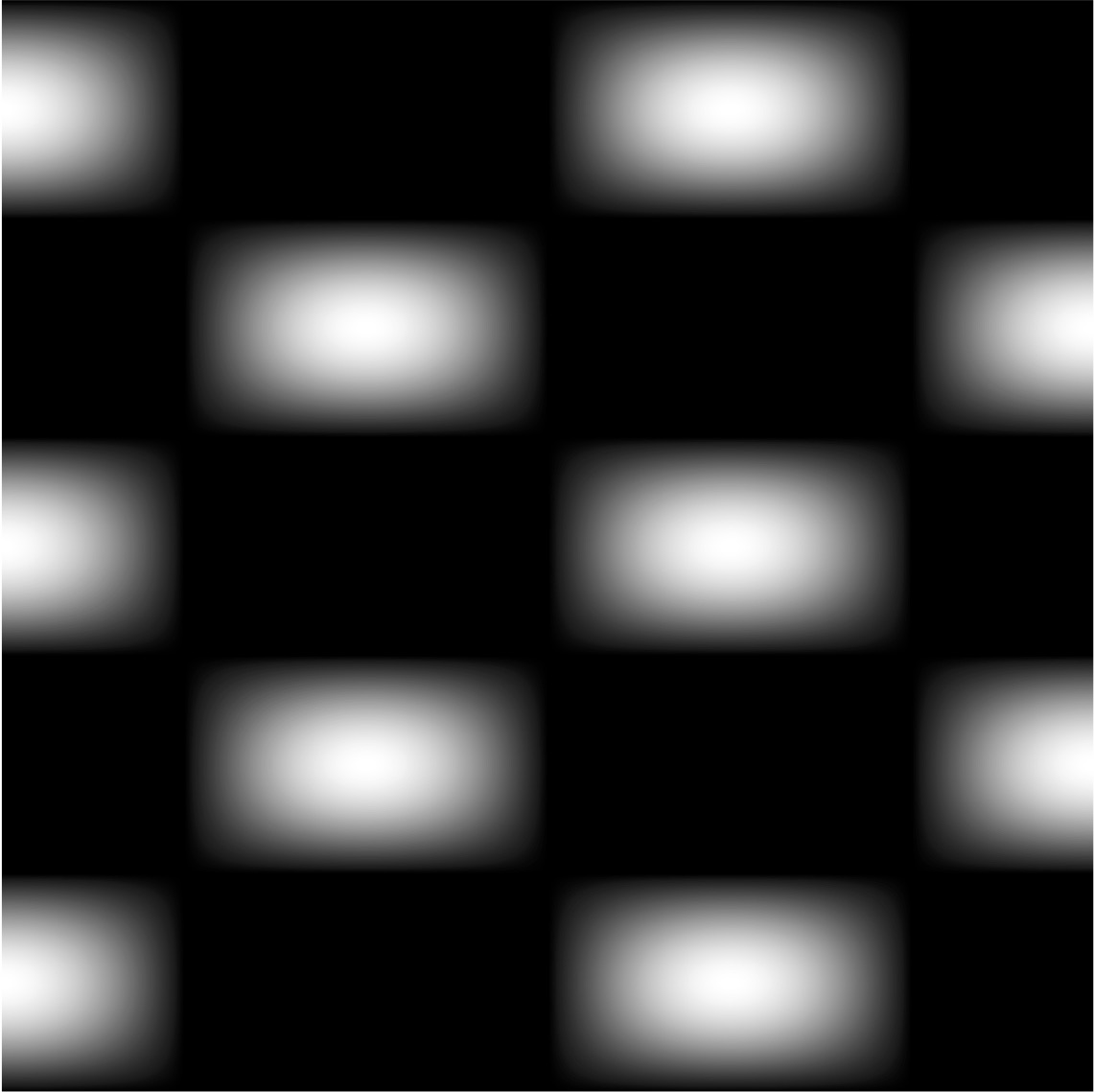}  
\end{center}
\caption{\small{$f(x,y)=\sin(5\pi x)\cos(3\pi y)\chi_{[0,1]^2}(x,y)$}}
\label{fig:cont_orig}
\end{figure}

\begin{figure}
\begin{tabular}{cl}
       & \hspace{1.3cm} TFS \hspace{1.9cm} GS with Haar \hspace{1.2cm} GS with DB2 \hspace{1.3cm} GS with DB3\\
{\rotatebox{90}{\hspace{1cm}SNR=0}}  &  \includegraphics[width=0.22\textwidth]{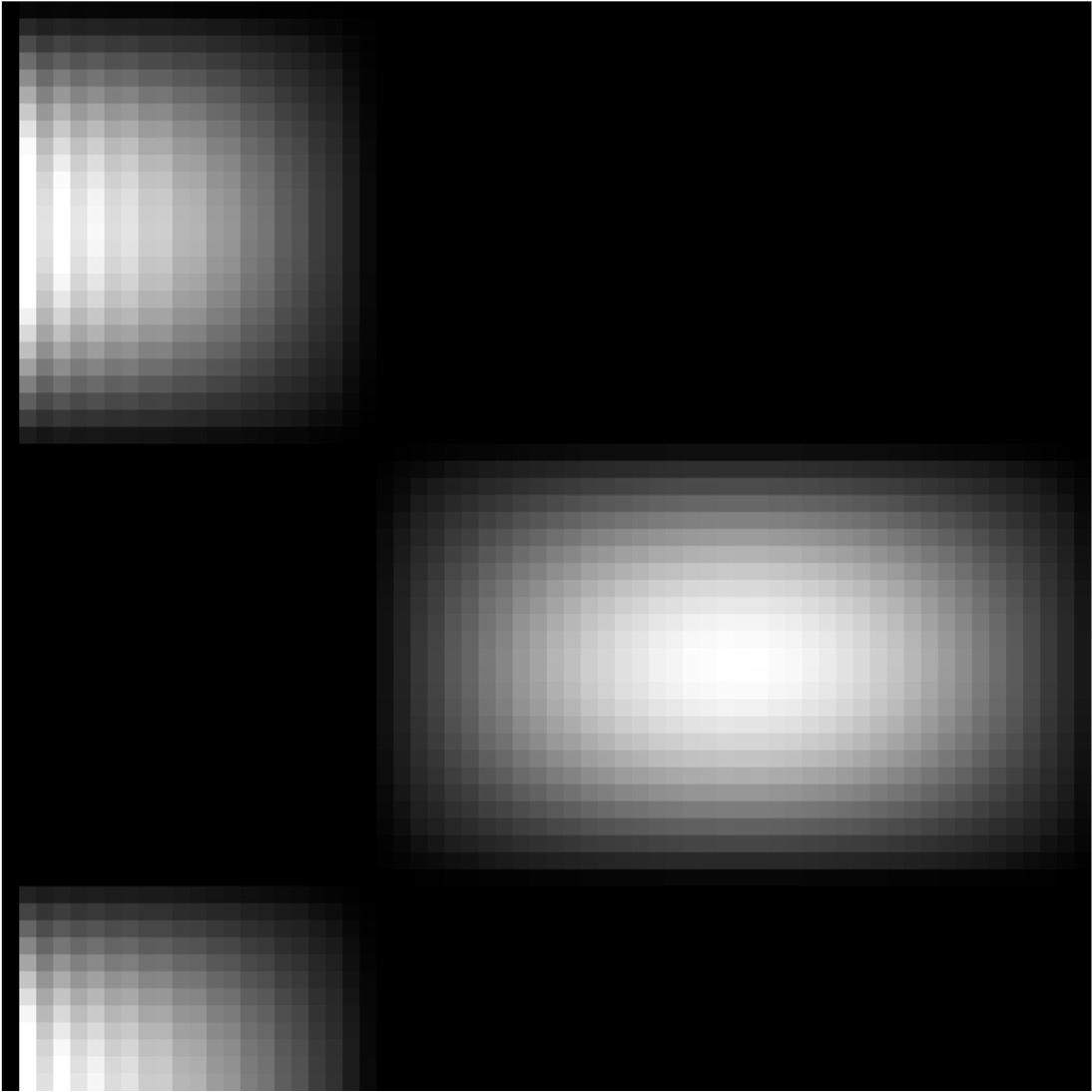}  \includegraphics[width=0.22\textwidth]{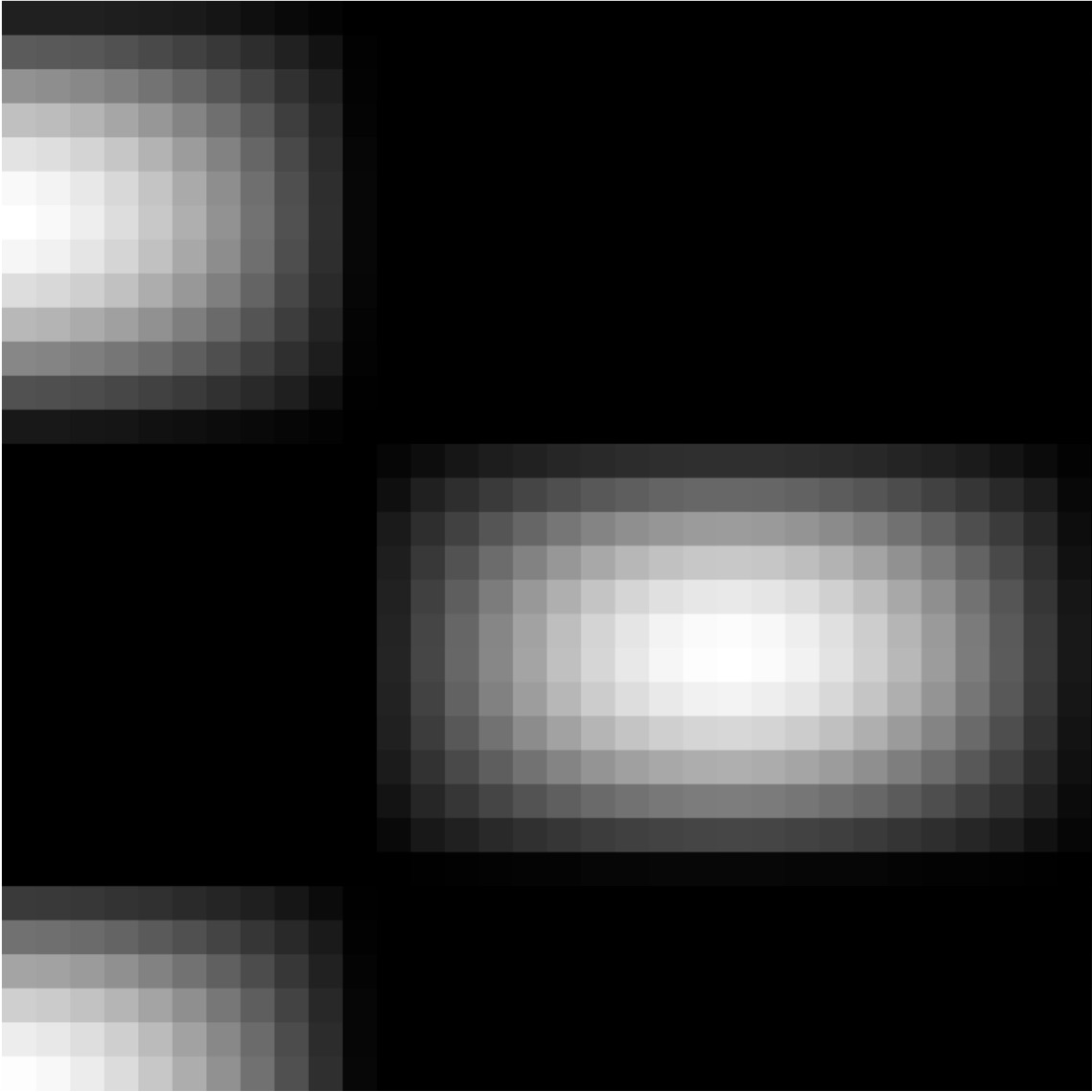}  \includegraphics[width=0.22\textwidth]{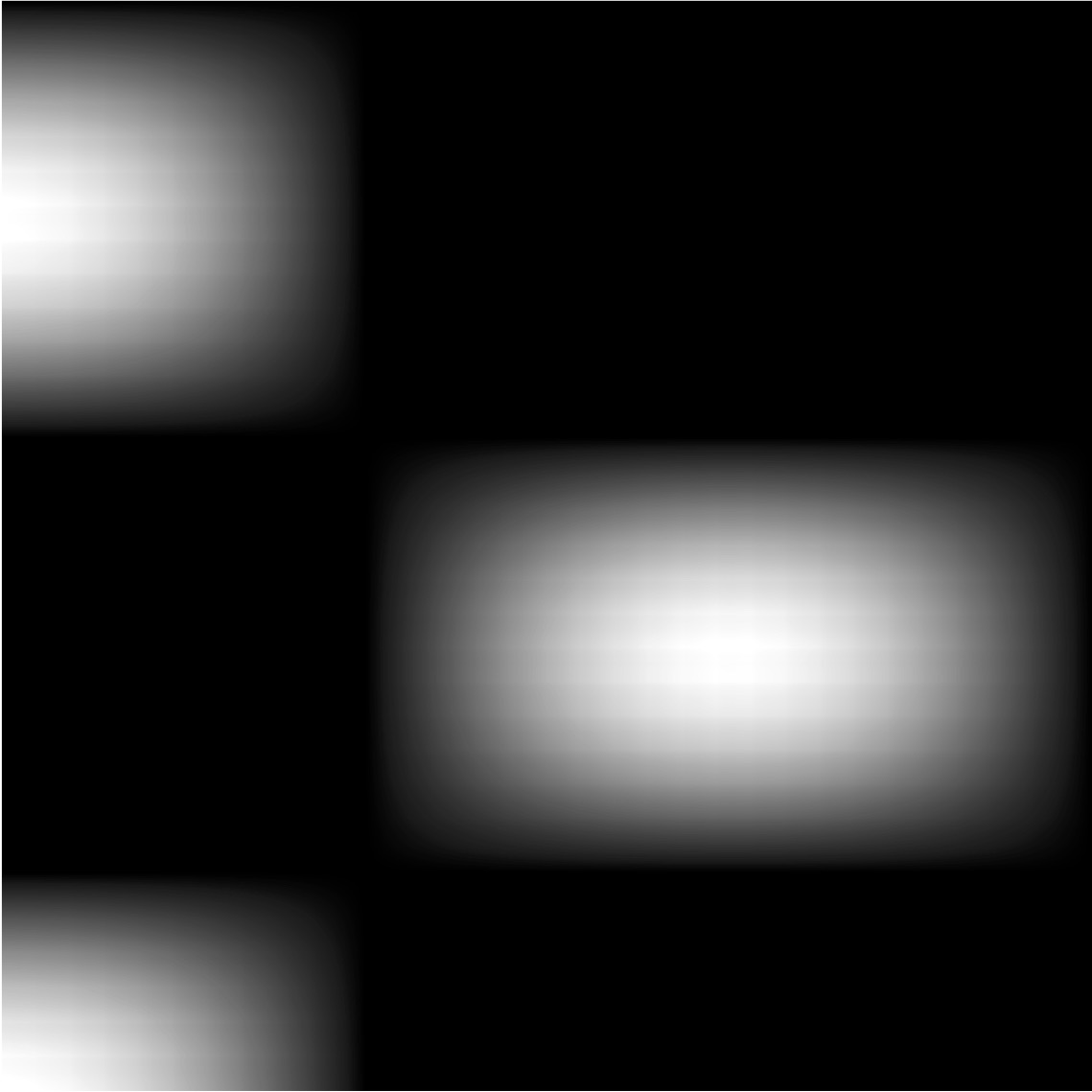}  \includegraphics[width=0.22\textwidth]{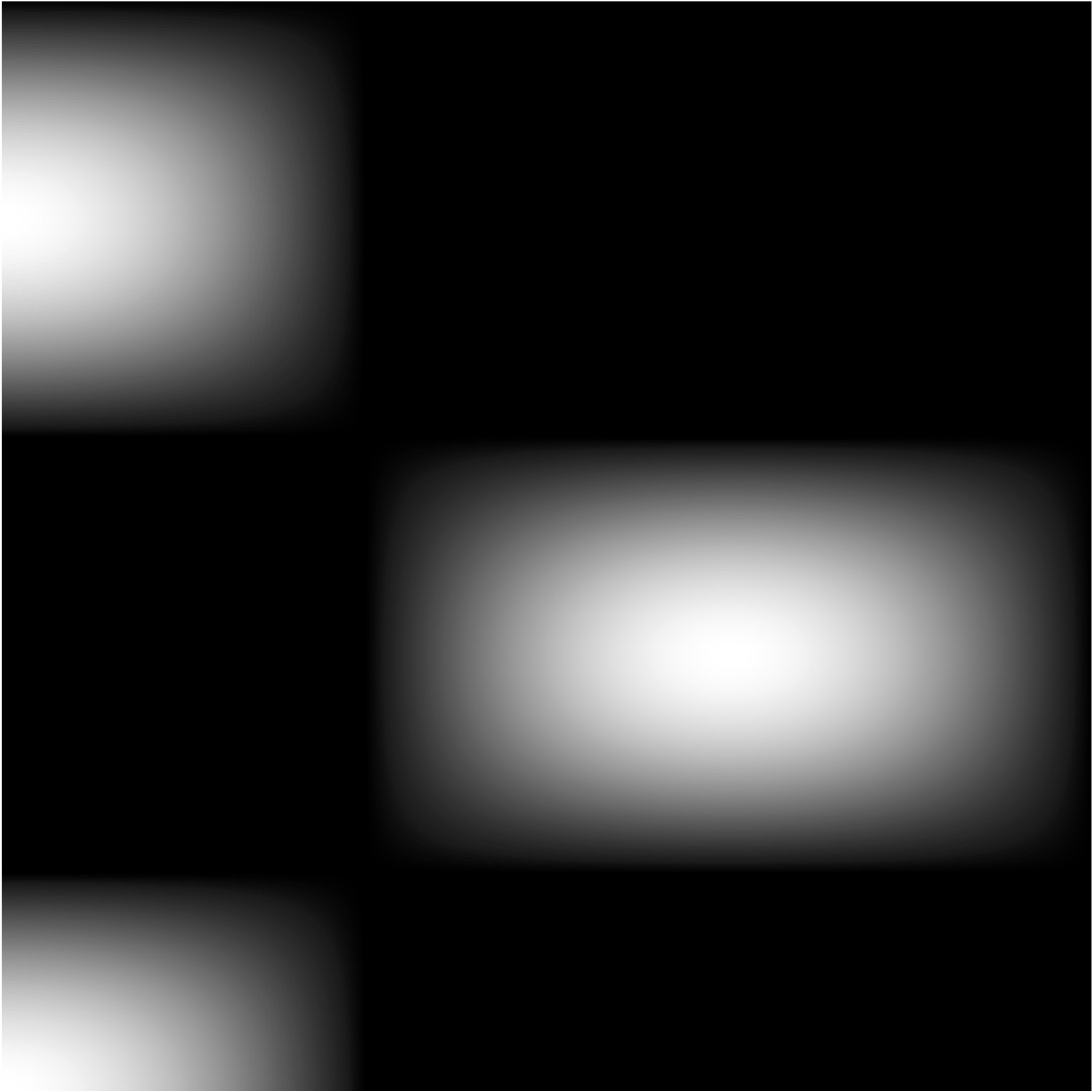}    \\
{\rotatebox{90}{\hspace{1cm}SNR=30}} &  \includegraphics[width=0.22\textwidth]{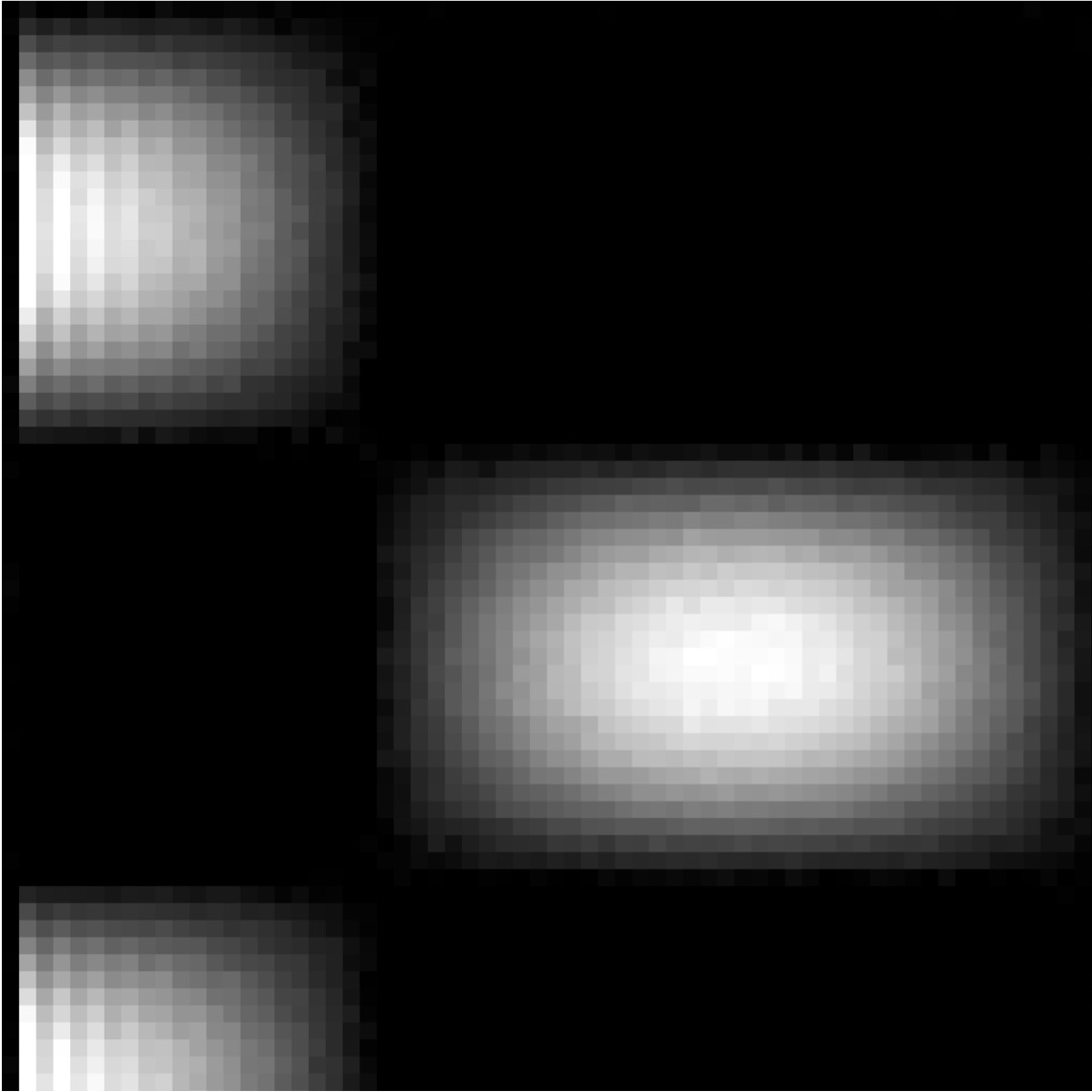} \includegraphics[width=0.22\textwidth]{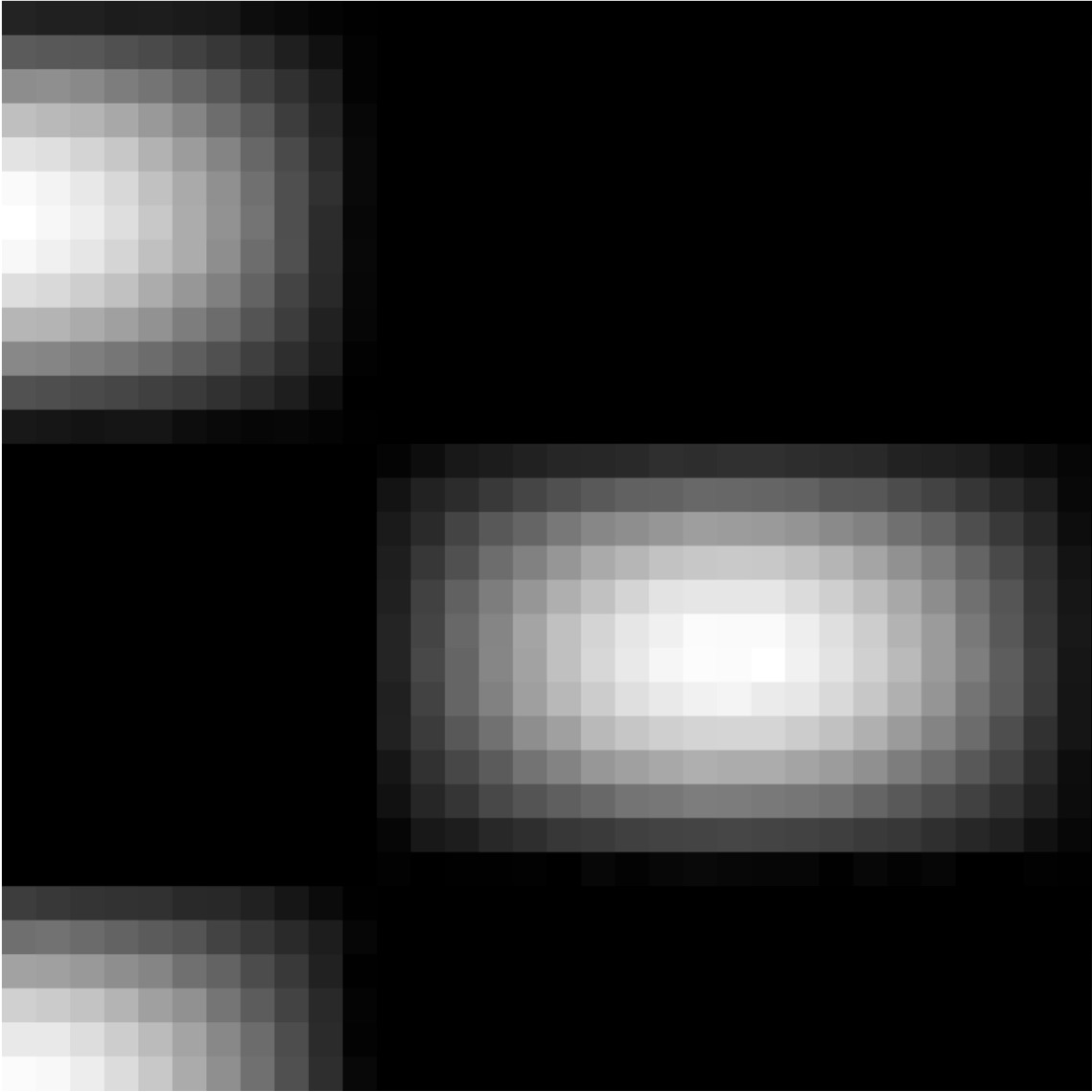}  \includegraphics[width=0.22\textwidth]{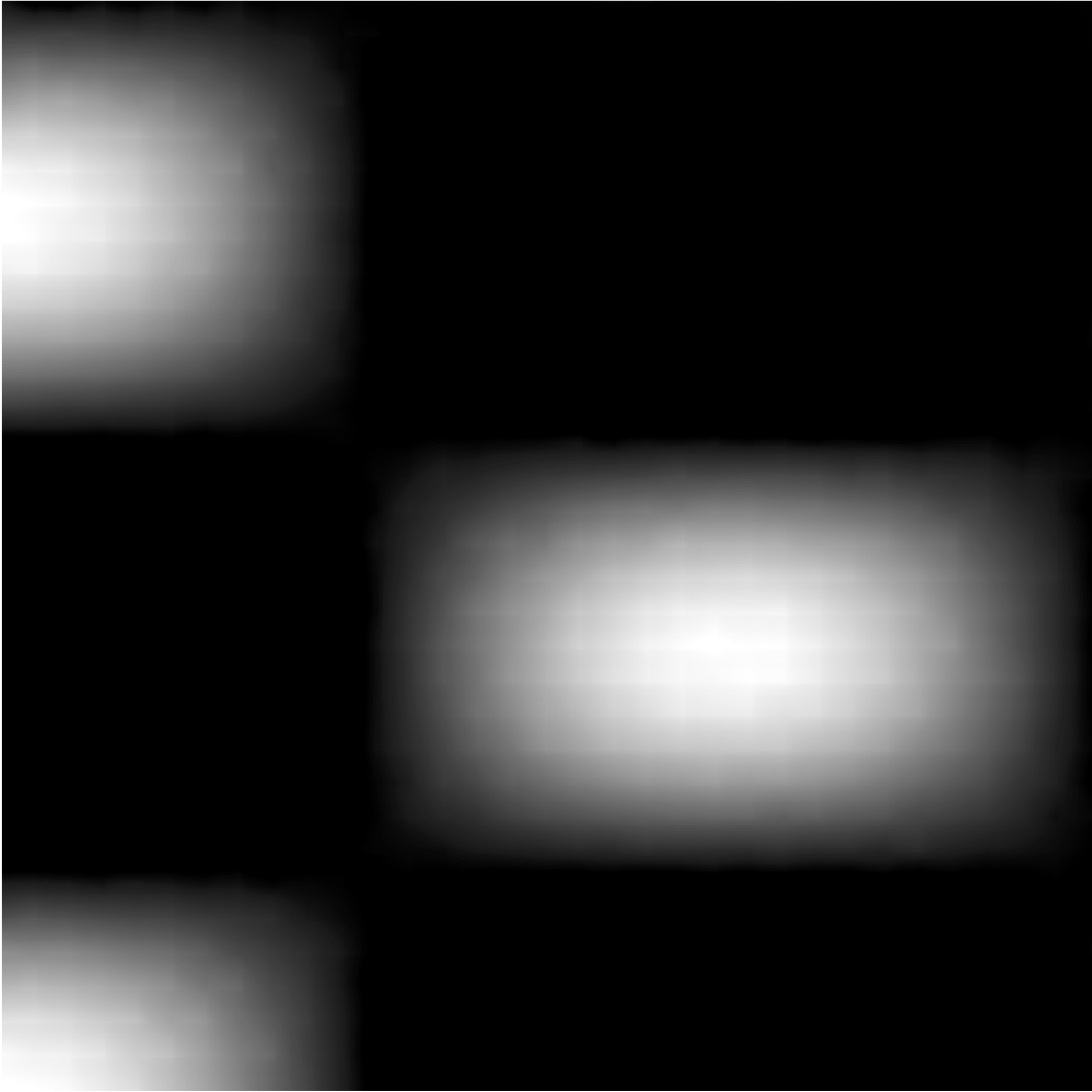}  \includegraphics[width=0.22\textwidth]{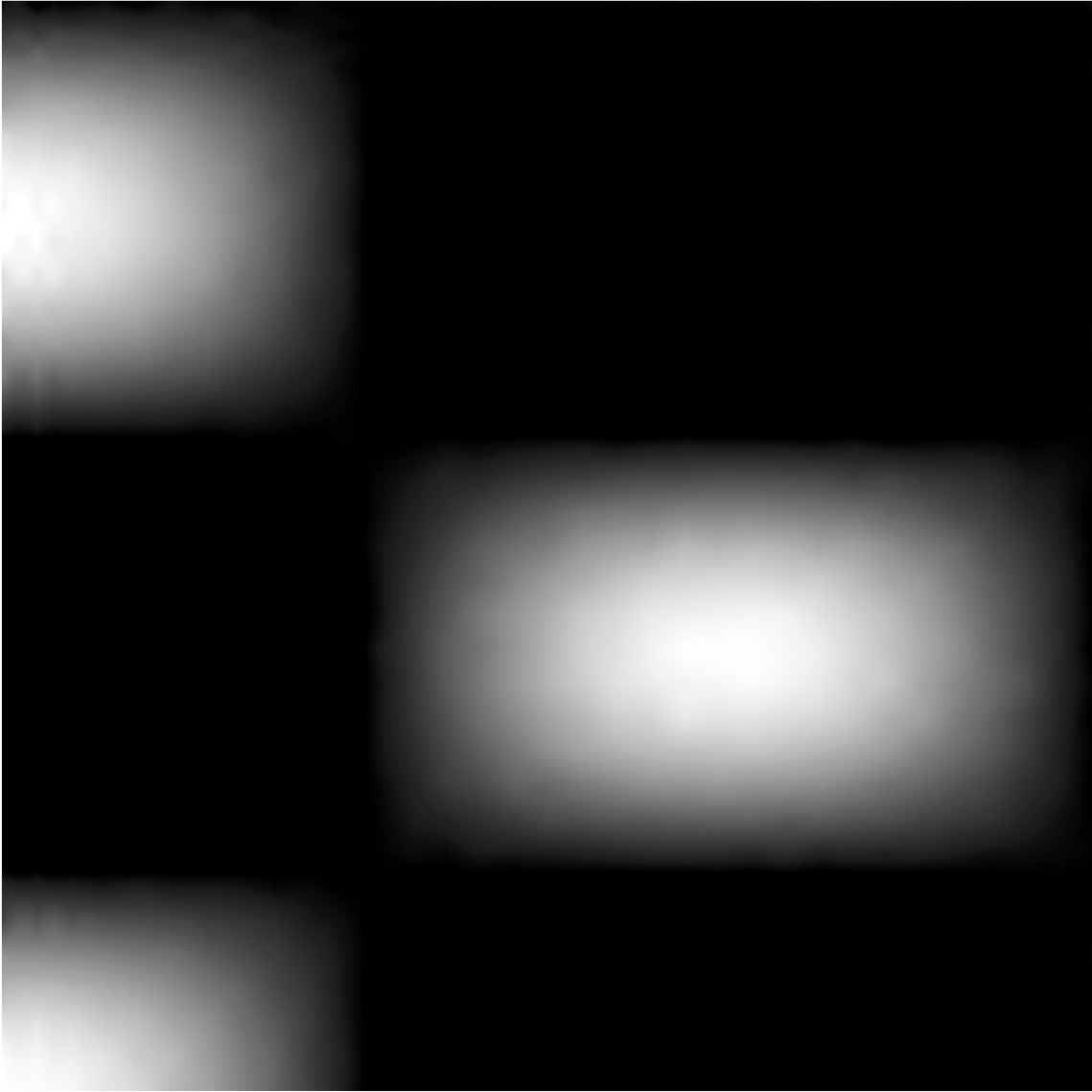} 
\end{tabular}
\caption{\small{
Top-left corner close-ups of the reconstructed function of $f$ in Figure \ref{fig:cont_orig} using uniform samples. 
}}
\label{fig:cont_uniform}
\end{figure}

\begin{figure}
\begin{tabular}{cl}
        & \hspace{1.5cm} Gr \hspace{1.9cm} GS with Haar \hspace{1.2cm} GS with DB2 \hspace{1.3cm} GS with DB3\\
{\rotatebox{90}{\hspace{1cm}SNR=0}}  & \includegraphics[width=.22\textwidth]{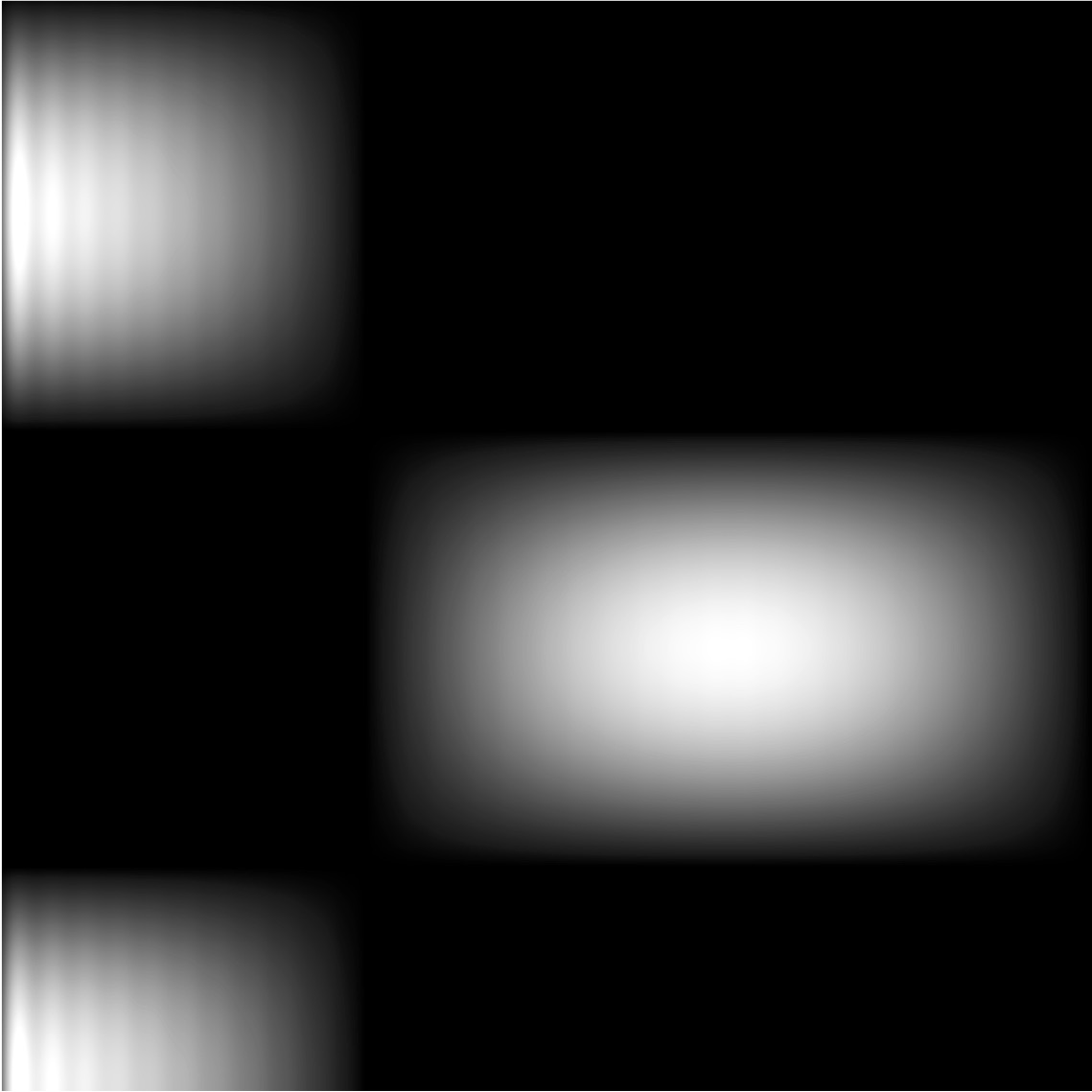}  \includegraphics[width=.22\textwidth]{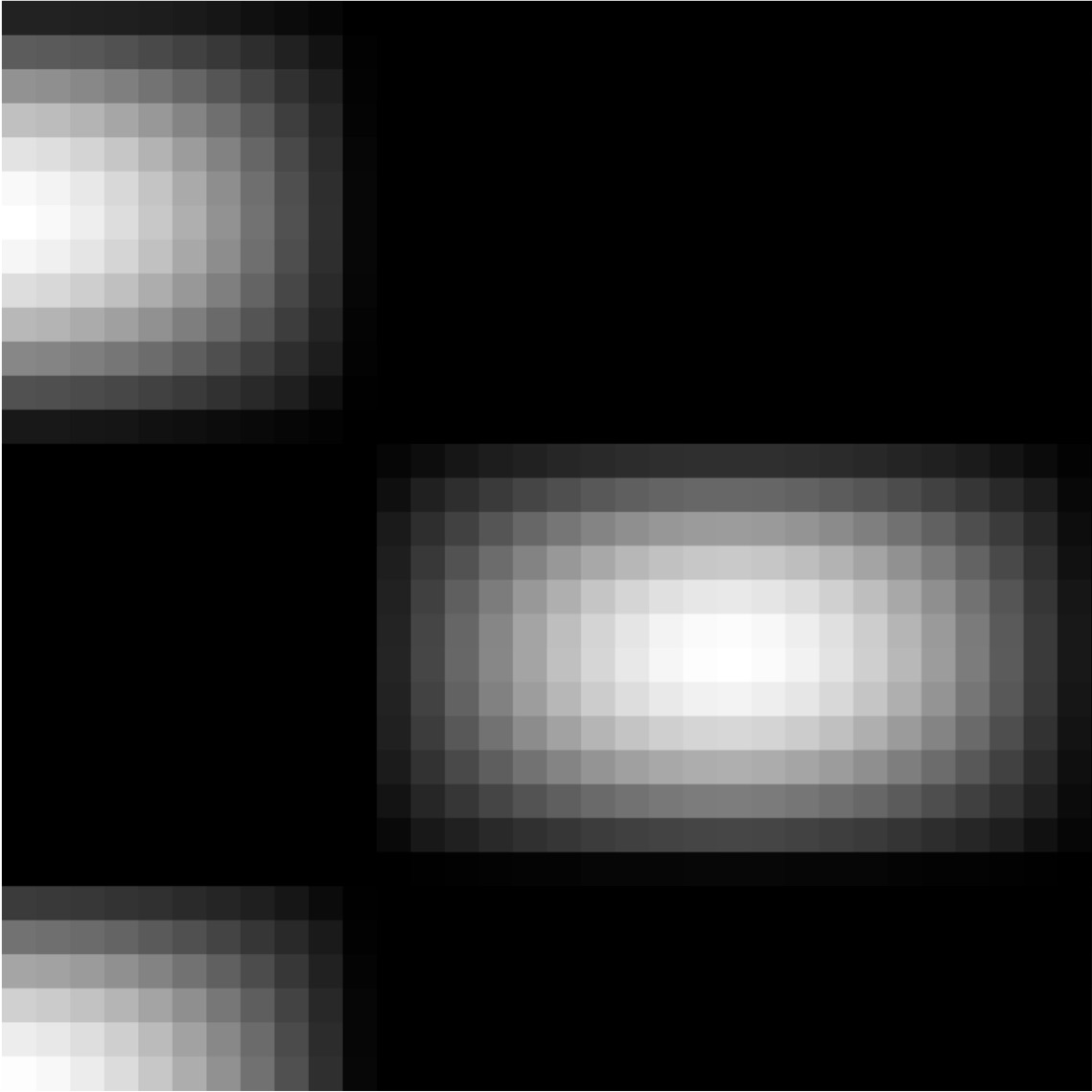}  \includegraphics[width=.22\textwidth]{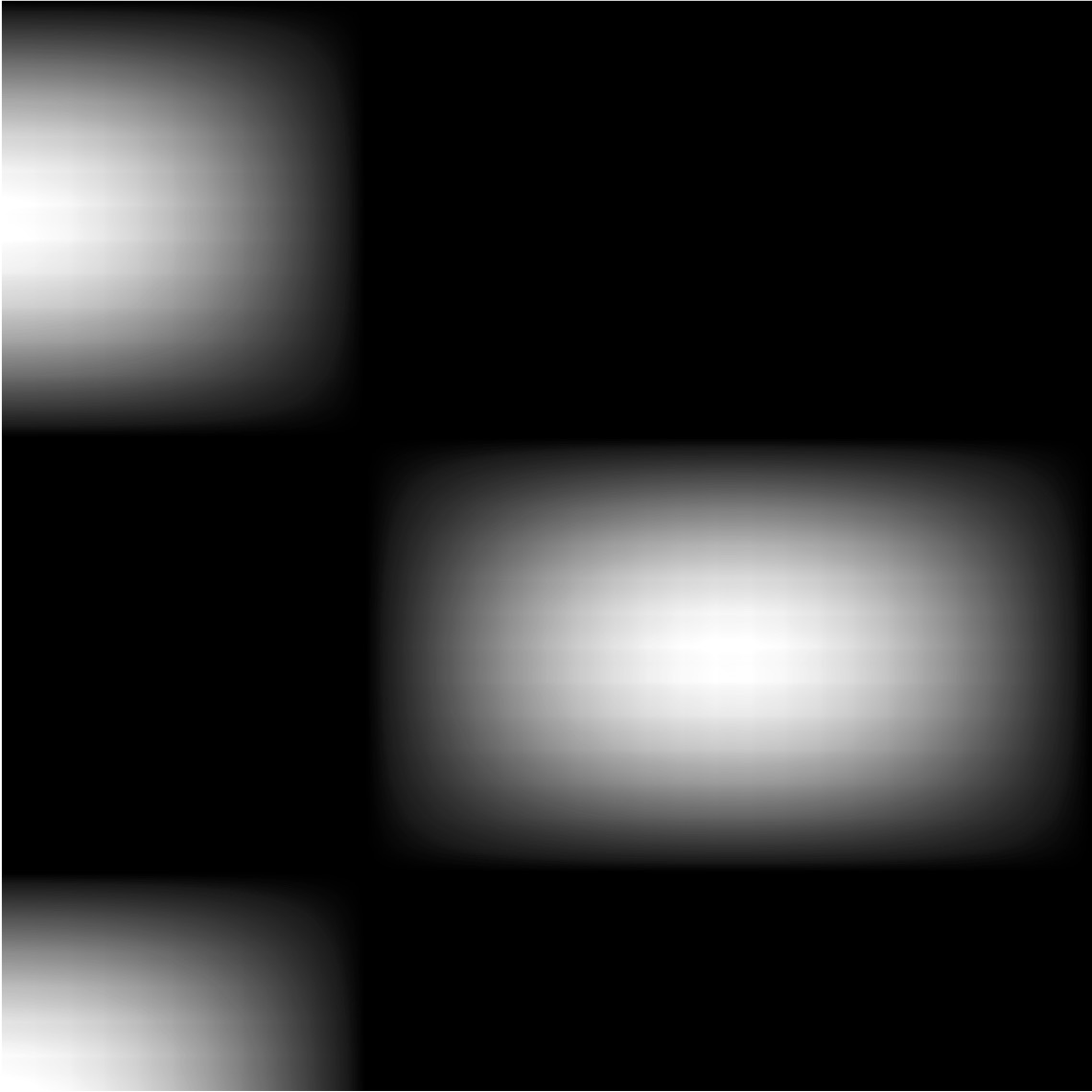}  \includegraphics[width=.22\textwidth]{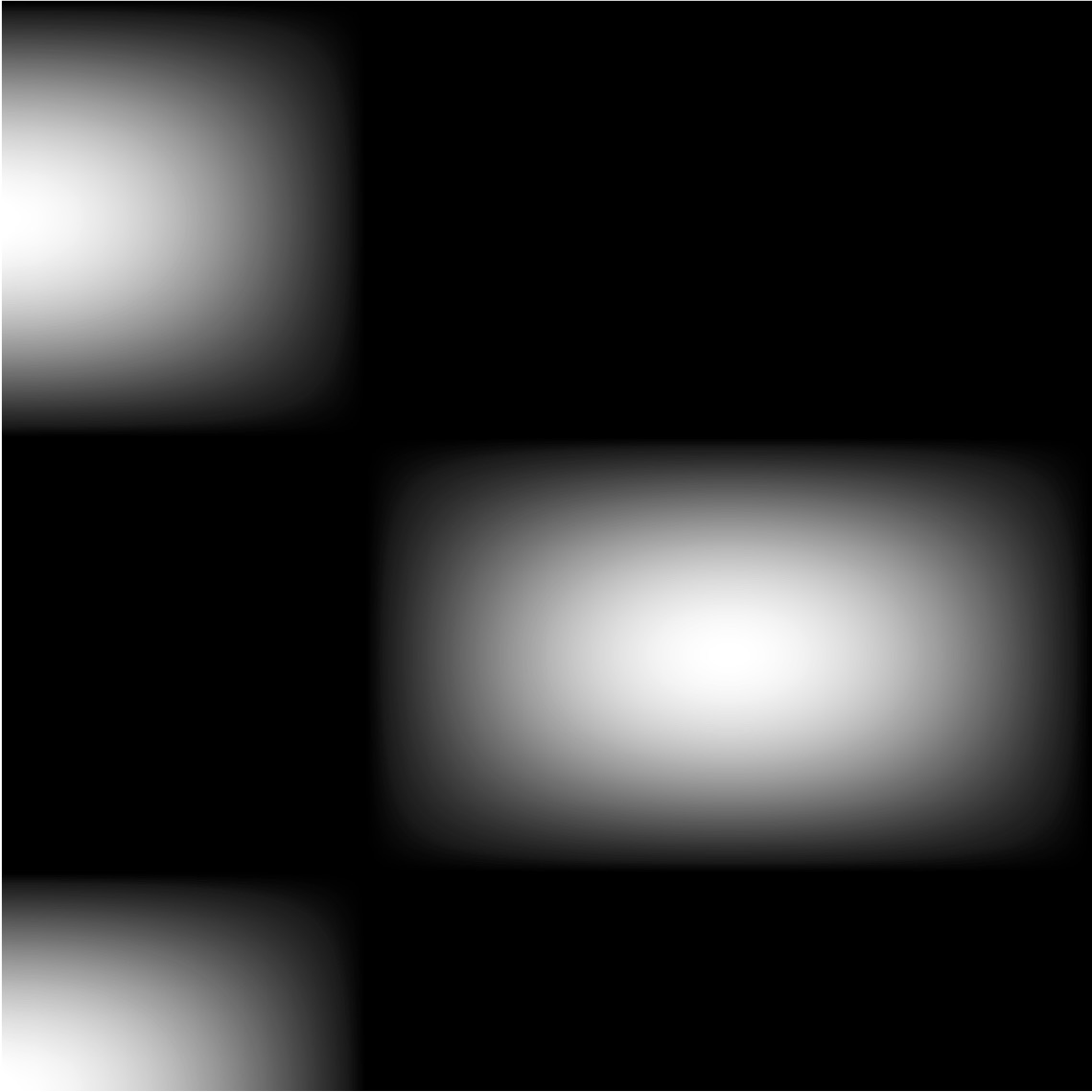}      \\
{\rotatebox{90}{\hspace{1cm}SNR=30}} & \includegraphics[width=.22\textwidth]{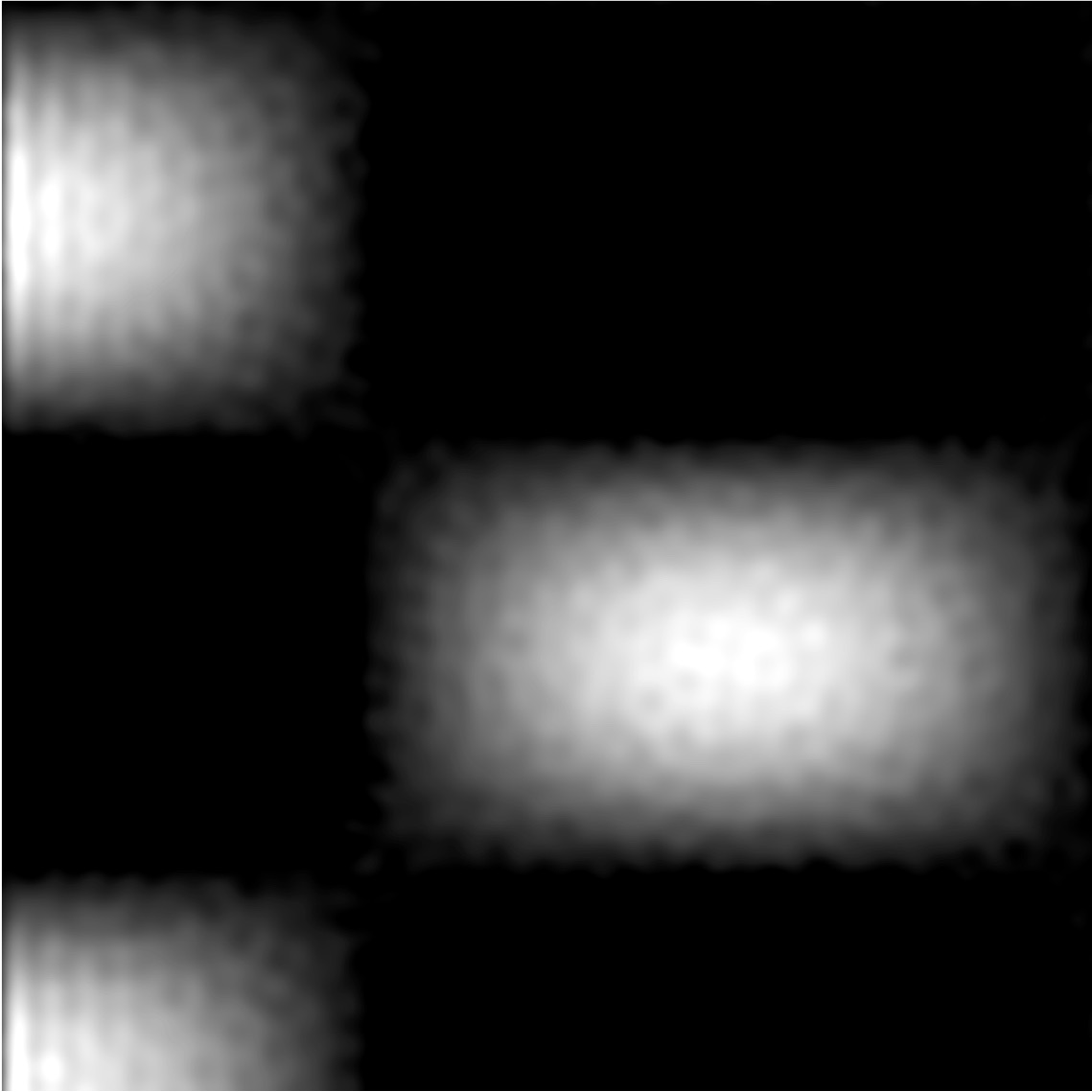} \includegraphics[width=.22\textwidth]{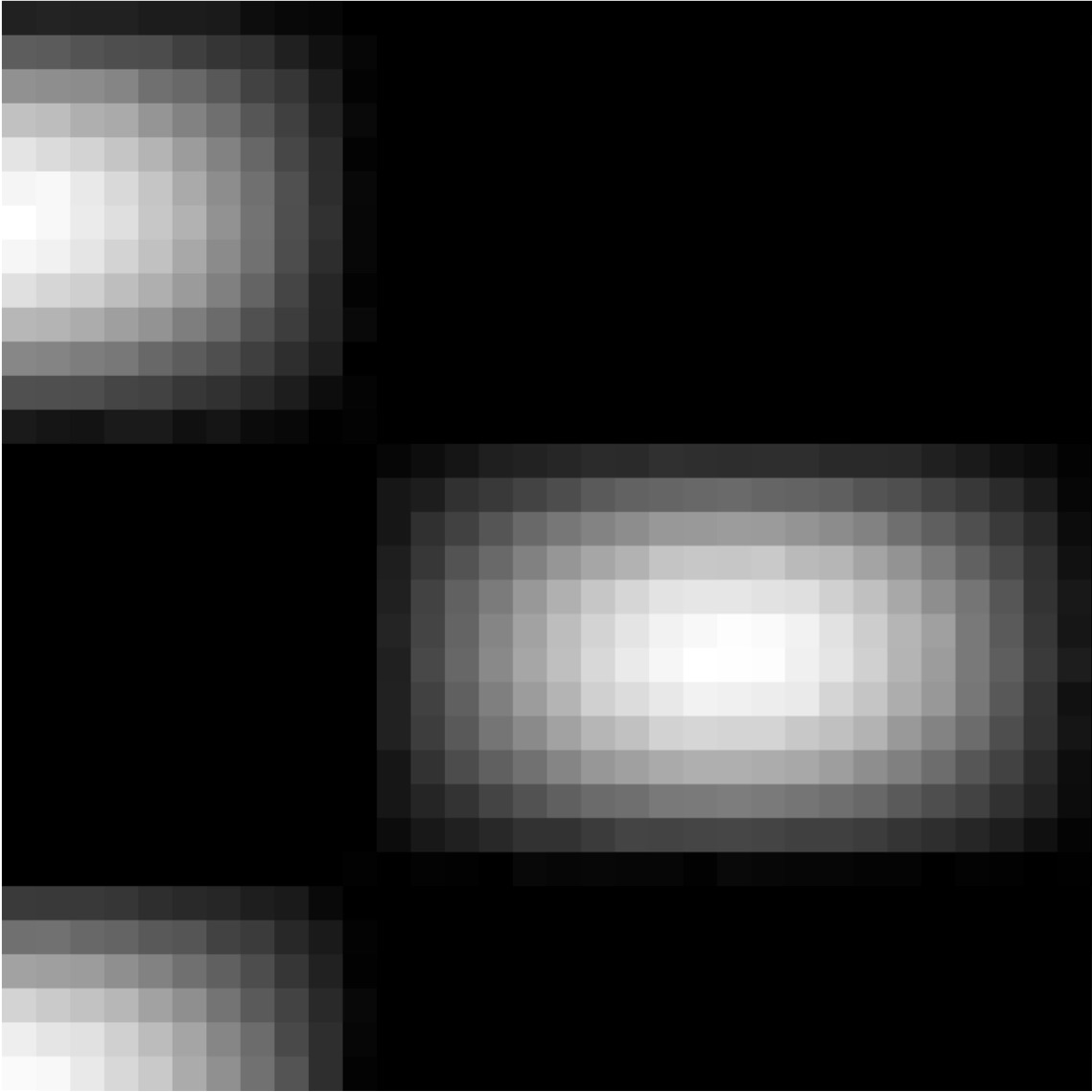}  \includegraphics[width=.22\textwidth]{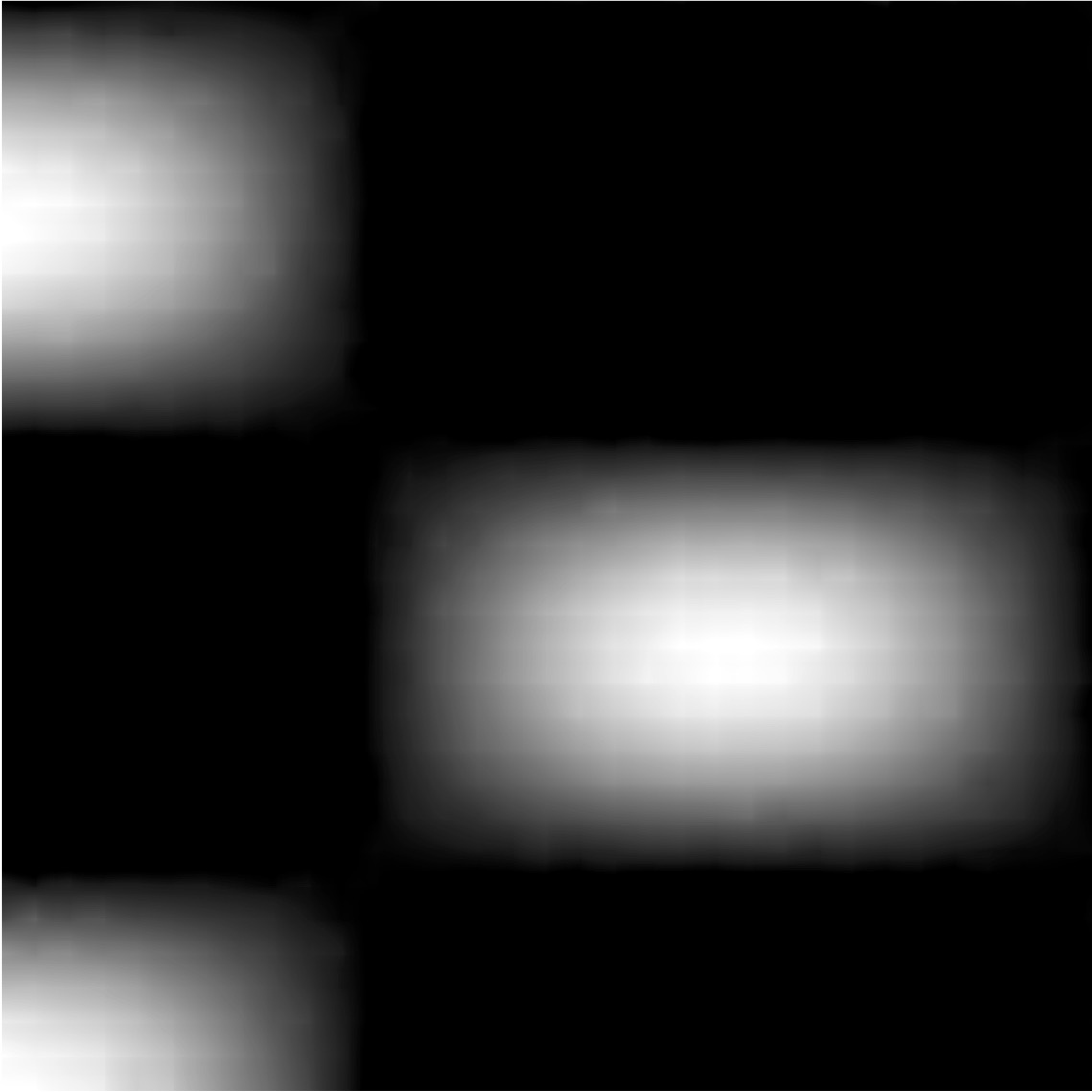}  \includegraphics[width=.22\textwidth]{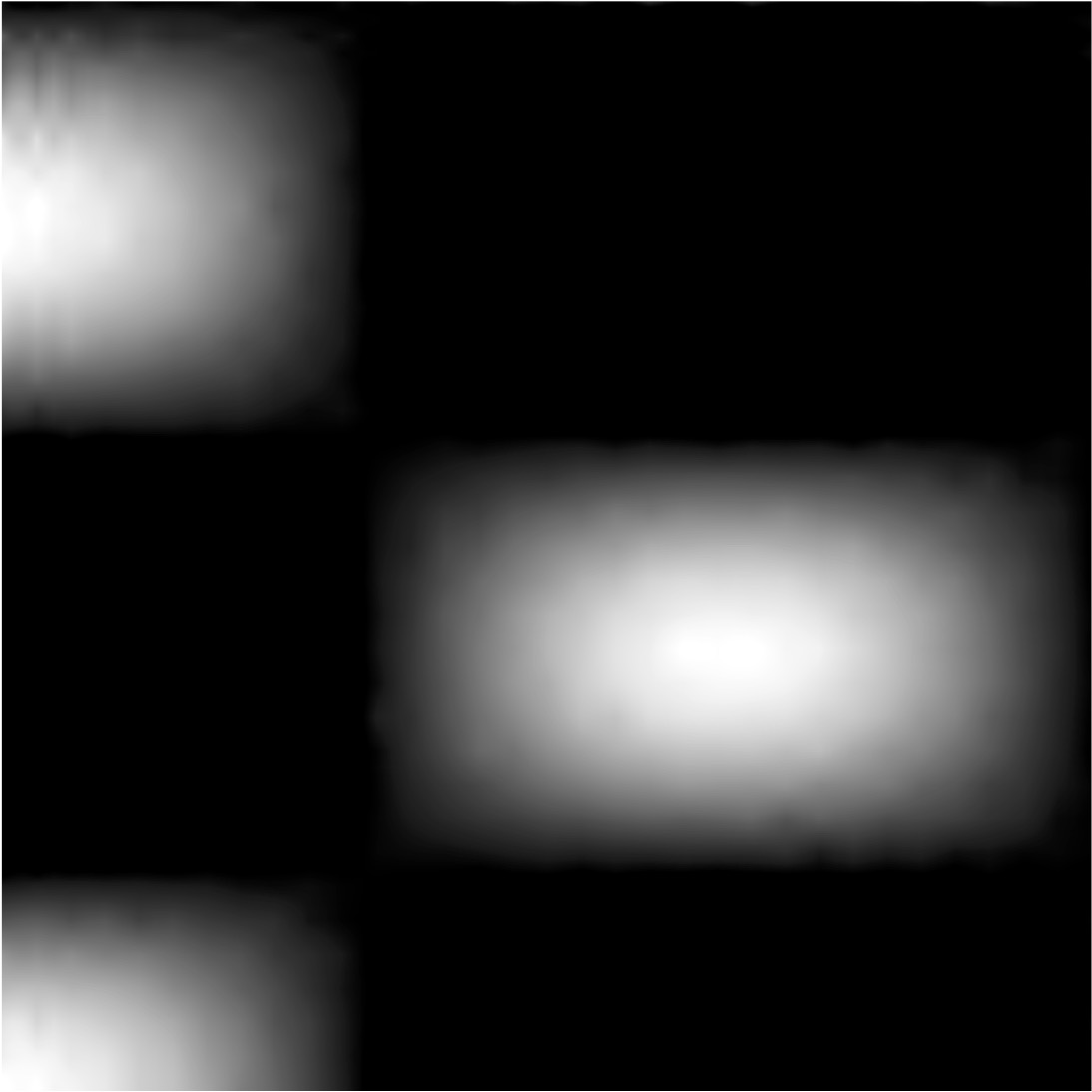}  
\end{tabular}
\caption{\small{
Top-left corner close-ups of the reconstructed function of $f$ in Figure \ref{fig:cont_orig} using radial samples. 
}}
\label{fig:cont_radial}
\end{figure}

\paragraph{Example 2.}
Next, we reconstruct a two-dimensional function shown in Figure \ref{fig:cont_orig}, which is again continuous but nonperiodic. We use an equispaced sampling scheme in Figure \ref{fig:cont_uniform}, and a radial sampling scheme that satisfies density condition \R{eq:density_con} in Figure \ref{fig:cont_radial}. By this example, we demonstrate robustness of GS when white Gaussian noise is added to the Fourier samples. We also demonstrate how one obtains improved reconstructions with increasing wavelet order.

\paragraph{Example 3.}
We report numerical error of reconstructions from Example 4 in Table \ref{t:error}. We see how error improves as the wavelet order increases.  As evident, the main advantage of the GS approach is the possibility of changing  the reconstruction space and using higher order wavelets for better performance.

\begin{table}
\begin{center}
\begin{tabular}{c|c|c|c|c|c|}
$\|f-\tilde f\|$         & SNR & TFS/Gr              & GS with Haar        & GS with DB2         & GS with DB3          \\ \hline
\multirow{2}{*}{Uniform} & 0   & $7.71\times10^{-2}$ & $4.13\times10^{-2}$ & $3.71\times10^{-3}$ & $8.11\times10^{-4}$  \\ \cline{2-6} 
                         & 30  & $7.88\times10^{-2}$ & $4.23\times10^{-2}$ & $9.14\times10^{-3}$ & $8.74\times10^{-3}$  \\ \hline
\multirow{2}{*}{Radial}  & 0   & $1.98\times10^{-2}$ & $4.13\times10^{-2}$ & $3.74\times10^{-3}$ & $7.95\times10^{-4}$  \\ \cline{2-6} 
                         & 30  & $2.93\times10^{-2}$ & $4.28\times10^{-2}$ & $1.07\times10^{-2}$ & $1.08\times10^{-2}$  \\ \hline
\end{tabular}
\caption{\small{
The $\cL_2$ error for the reconstruction of function $f(x,y)=\sin(5\pi x)\cos(3\pi y)\chi_{[0,1]^2}(x,y)$ via truncated Fourier series (TFS) or gridding (Gr), and via GS with $64\times64$ Haar, DB2 or DB3 boundary-corrected wavelets from noiseless (SNR=0) and noisy (SNR=30) pointwise Fourier measurements taken on a uniform or radial sampling scheme from the region $[-64,64]^2$. See Figures \ref{fig:cont_orig},\ref{fig:cont_uniform} and \ref{fig:cont_radial}.
}\label{t:error}}
\end{center}
\end{table}

\paragraph{Example 4.}

\begin{figure}
\begin{center}
\begin{flushleft}
\hspace{2.6cm} Original \hspace{3.1cm} TFS \hspace{2.7cm} GS with DB4 \end{flushleft}
\includegraphics[width=0.3\textwidth]{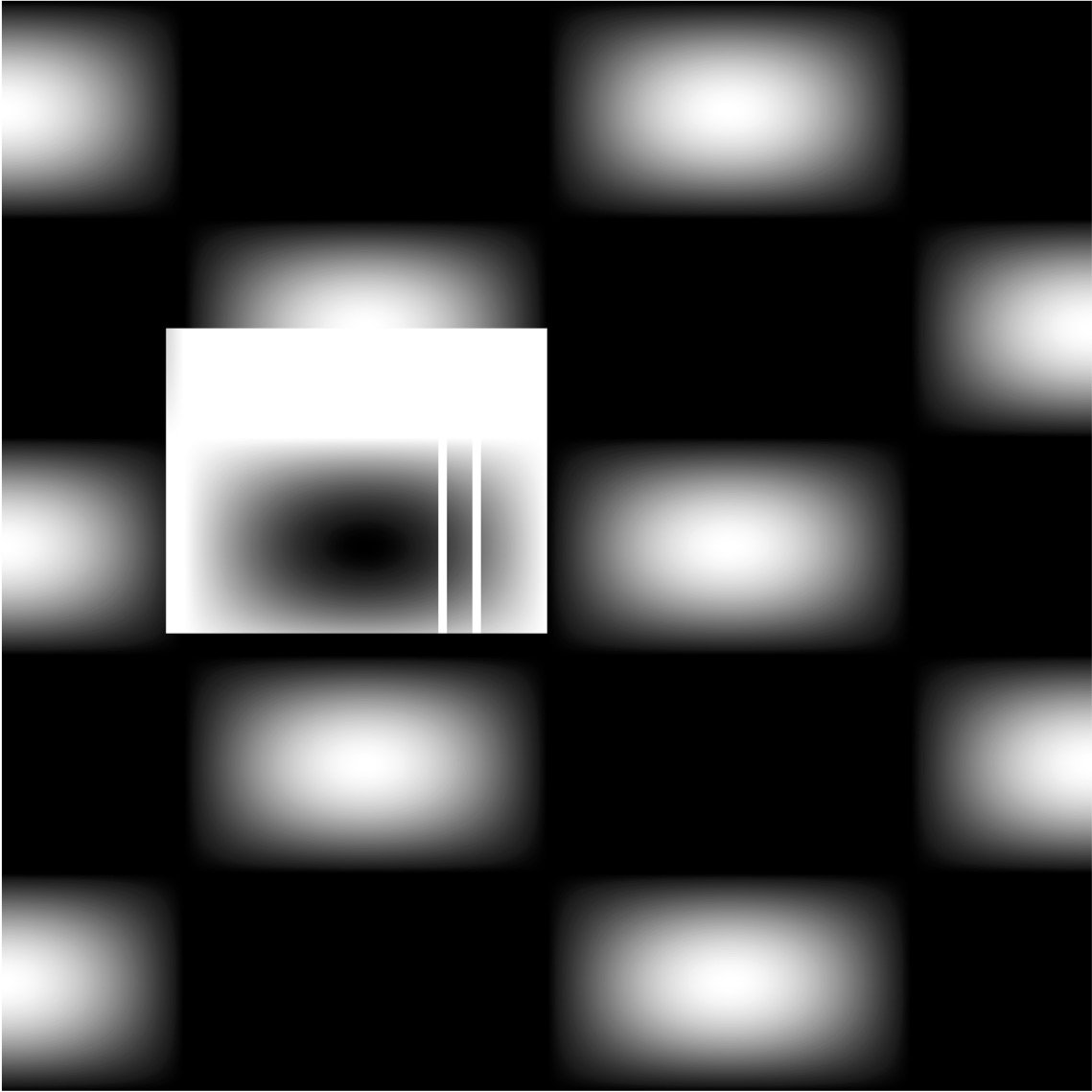}  \includegraphics[width=0.3\textwidth]{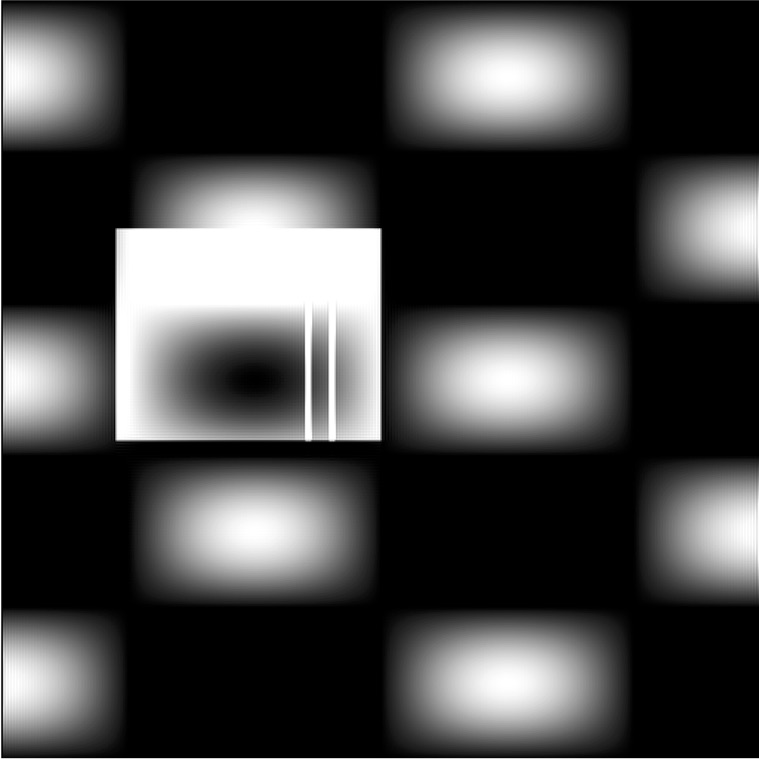} \includegraphics[width=0.3\textwidth]{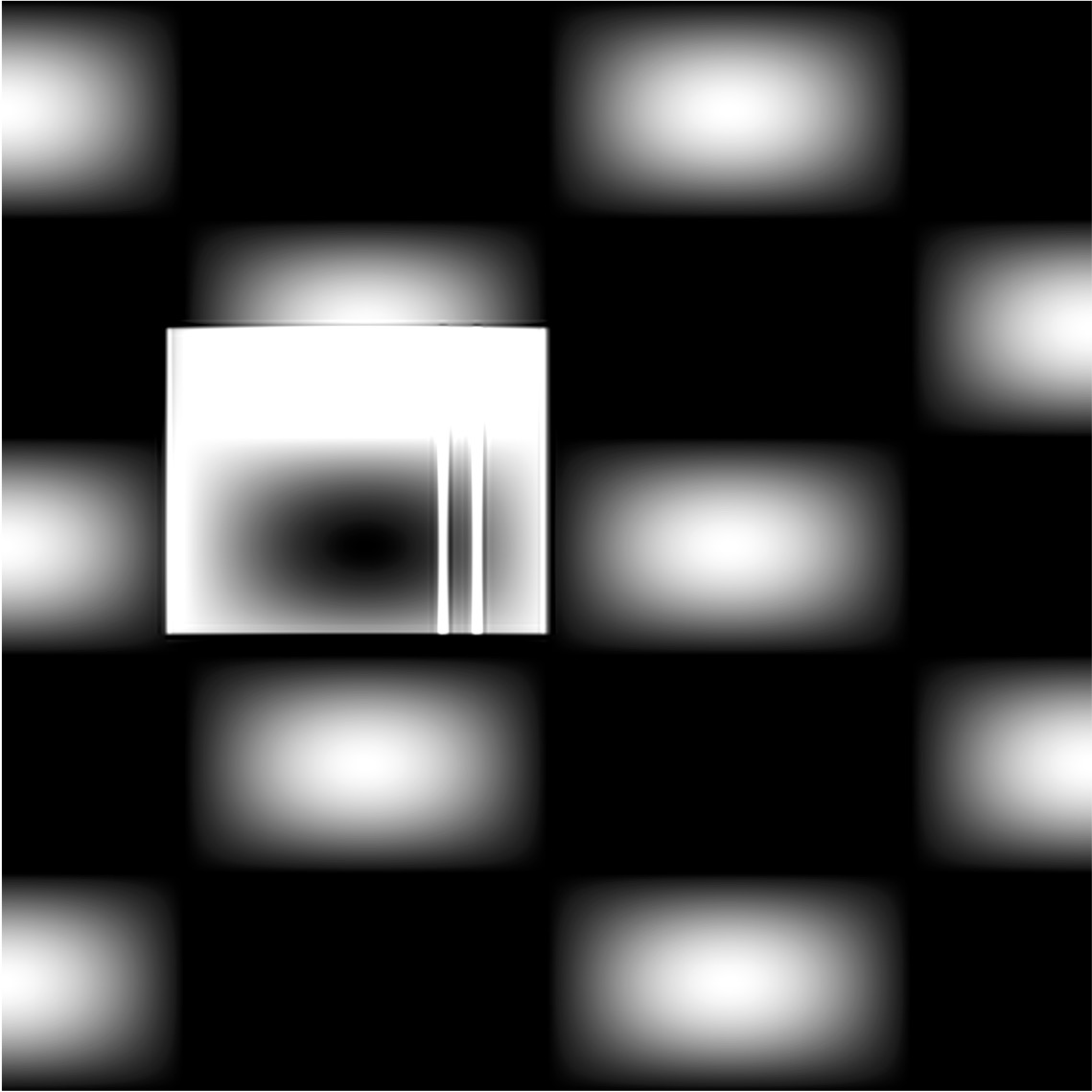}  \\
\includegraphics[width=0.3\textwidth]{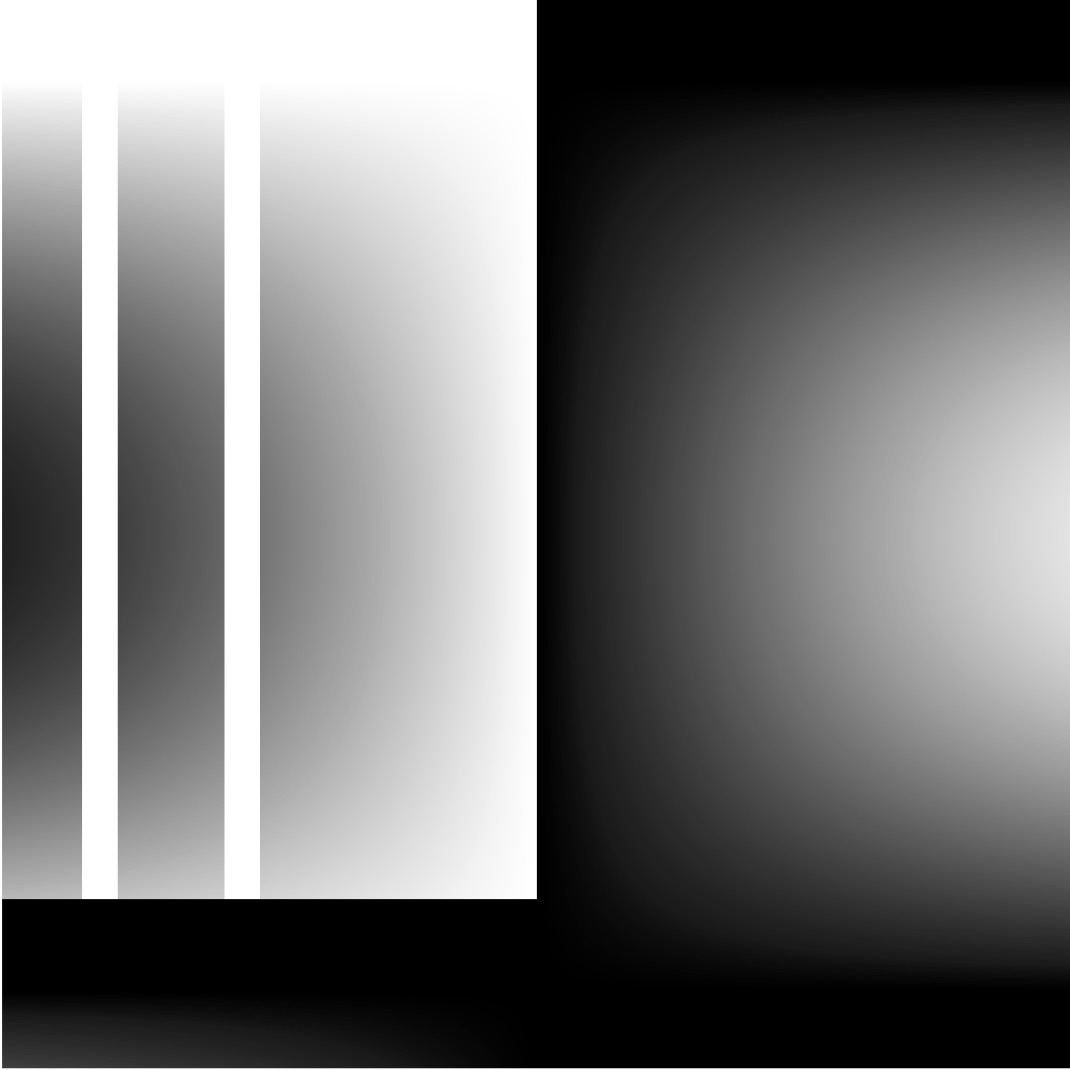}   \includegraphics[width=0.3\textwidth]{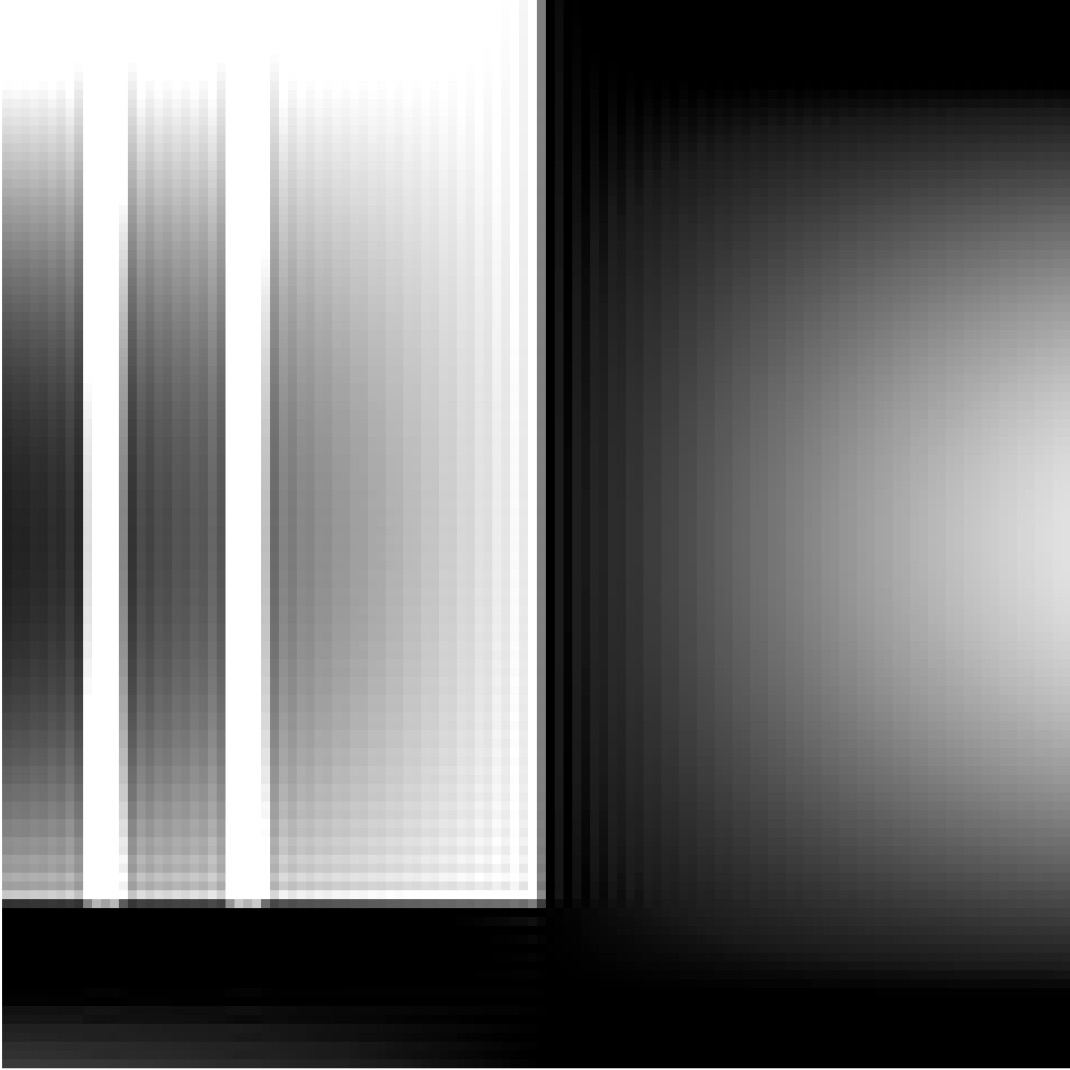} \includegraphics[width=0.3\textwidth]{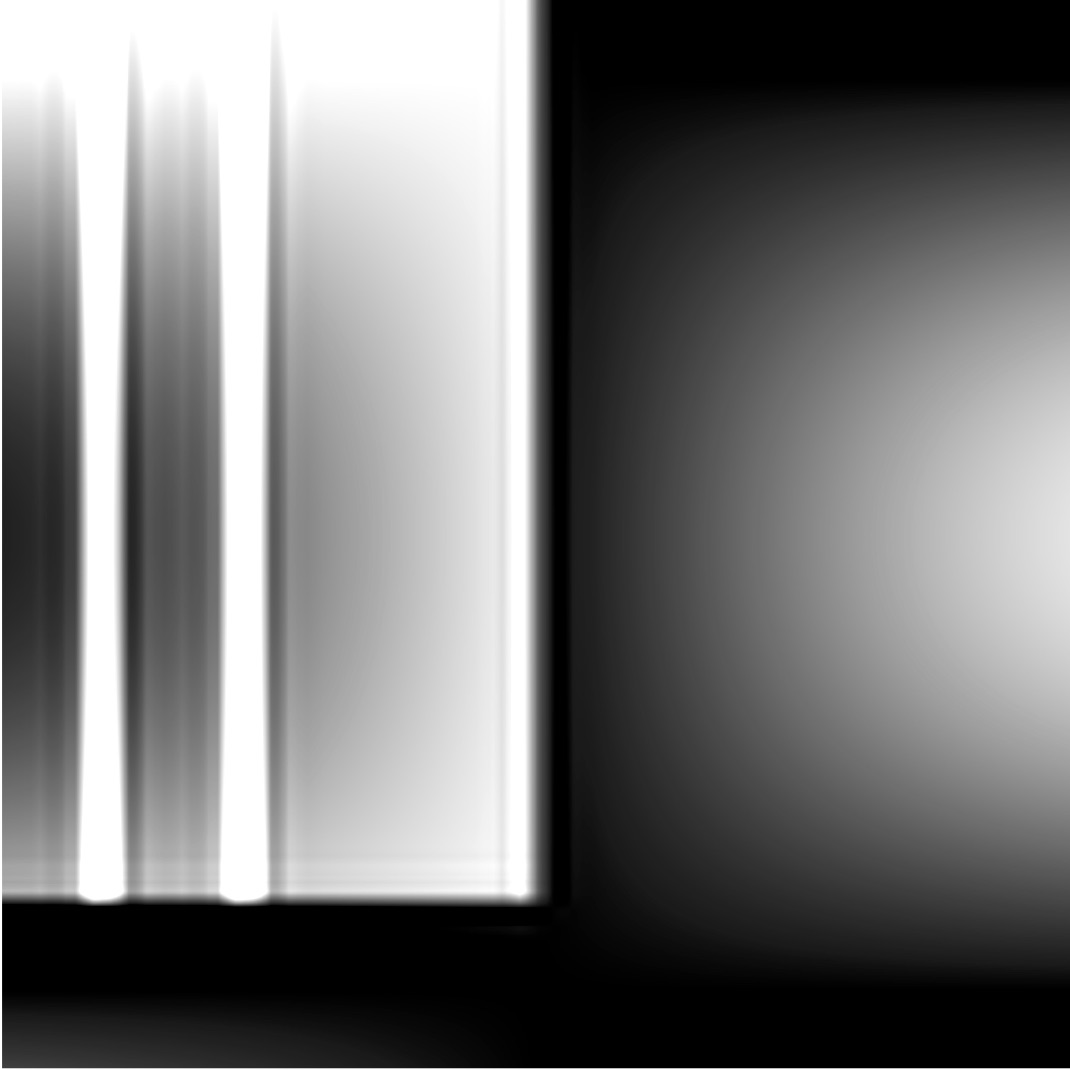}
\end{center}
\caption{\small{A discontinuous function is reconstructed by truncated Fourier series and by GS using $256\times256$ DB4 boundary-corrected wavelets from $512\times512$ pointwise Fourier measurements taken from on an equispaced grid.
}}
\label{fig:dicont}
\end{figure}

In Figure \ref{fig:dicont}, we demonstrate reconstruction of a discontinuous function via GS with DB4, and compare it to the direct approach, when using uniform samples.

\begin{figure}
\begin{center}
\begin{tabular}{@{\hspace{0pt}}c@{\hspace{-15pt}}c@{\hspace{0pt}}}
 Original& Sampling mask\\
\includegraphics[width = 0.4\textwidth]{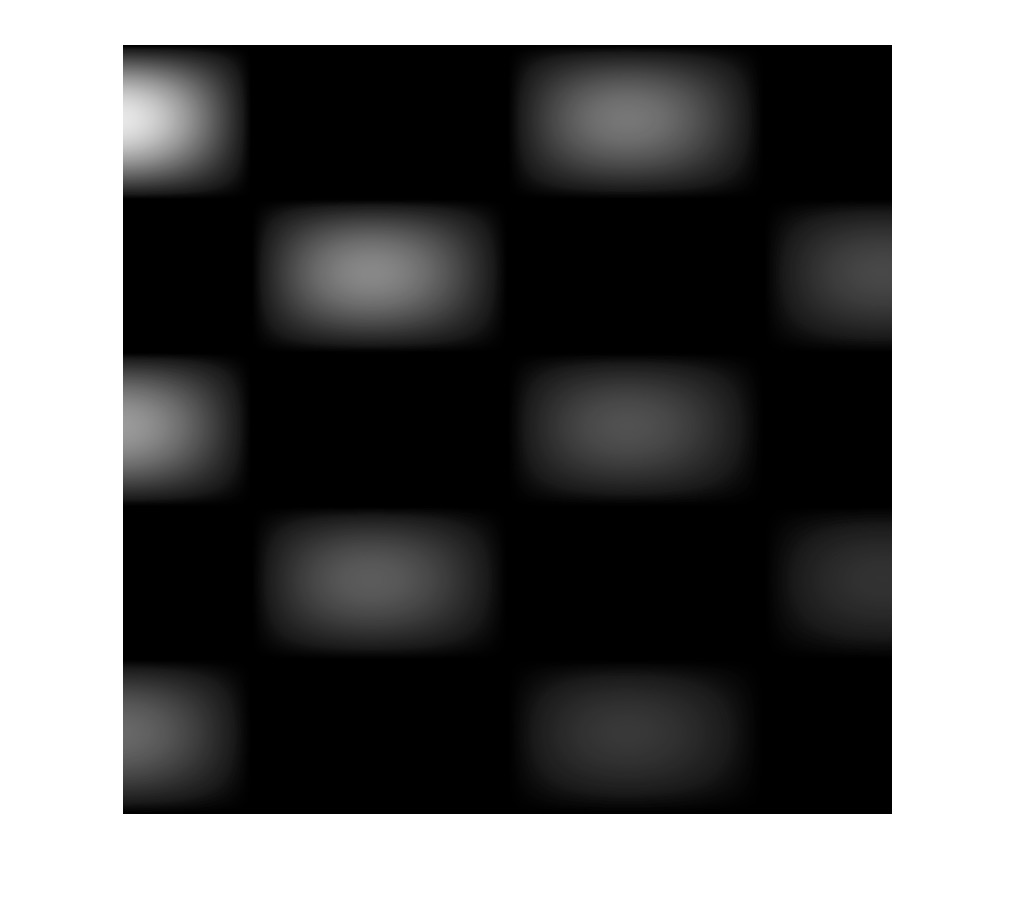} &\includegraphics[width = 0.4\textwidth]{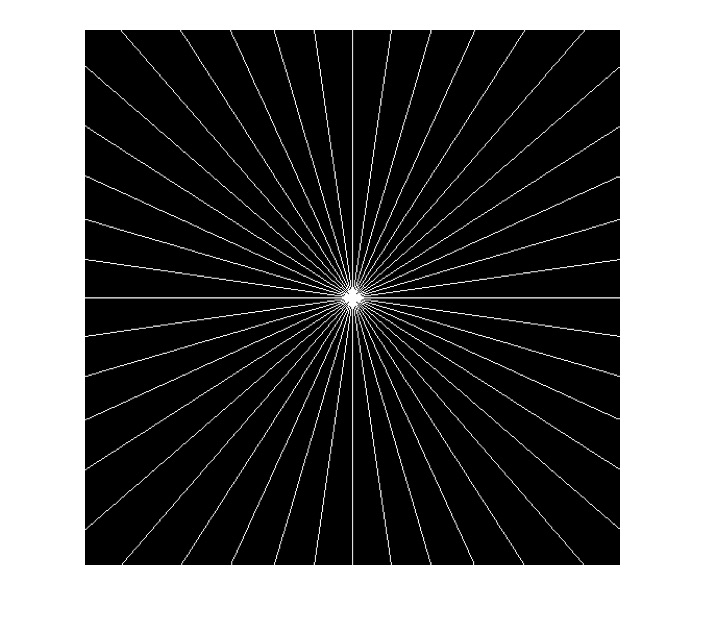} \\
GS with $\ell^1$ & CS with DFT and DWT\\
\includegraphics[width = 0.47\textwidth]{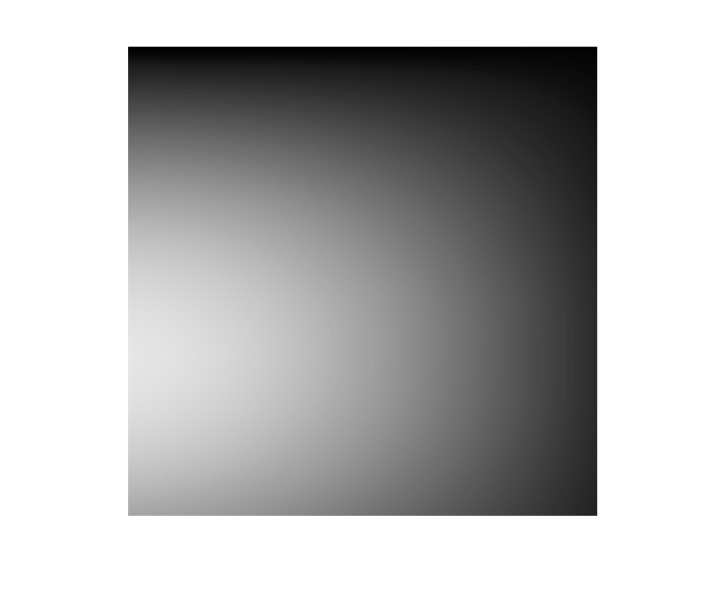} & \includegraphics[width = 0.47\textwidth]{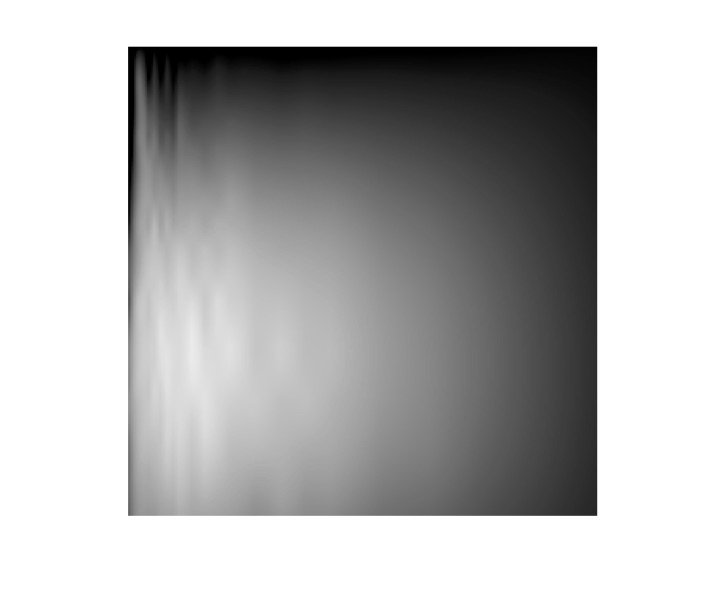}
\end{tabular}
\end{center}
\caption{\small{
The reconstruction of a smooth function (top left) from $4.25\%$ of its first $512\times 512$ Fourier coefficients using a standard compressed sensing (CS) approach with DFT and DWT, and using the $\ell^1$ generalized sampling (GS) approach. The samples taken for both reconstructions are those supported on the start shape sampling mask (top right). The bottom row shows a zoom-in of the top left reconstruction. The error from the standard CS approach is $1.6\times 10^{ -2}$, while the error from the GS approach is $4.7\times 10^{-3}$.
} \label{fig:CS}
}
\end{figure}

\paragraph{Example 5.}
Figure \ref{fig:CS} illustrates the advantage of combining generalized sampling with $\ell^1$ minimization from \R{eq:gsl1_fin}, as opposed to simply employing standard finite dimensional compressed sensing techniques where one assumes compressibility with respect to some wavelet bases and one solves
\be{\label{eq:cs_fin}
x \in \mathrm{arg} \min_{z\in\bbC^N} \nm{z}_{\ell^1} \text{ subject to } \nm{P_\Omega U_{df} U_{dw}^{-1} z - y}_{\ell^1}\leq \delta,
}
for given (noisy) Fourier samples $y$. The reconstruction is then taken to be $U_{dw}^{-1} x$. Here, $U_{df}$ is the discrete Fourier transform, $U_{dw}$ is the discrete wavelet transform and $P_\Omega$ is a projection matrix setting entries which are not indexed by $\Omega$ to zero. However, as explained in \cite{BAACHGSCS}, any solution to this will be limited by the error of the truncated Fourier representation of $f$, whereas the error of the solution to \R{eq:gsl1_fin} can be controlled by the decay of the underlying wavelet coefficients.

In this example, we consider the recovery of the smooth function
$$
f(x) = \cos(3x)\sin(5y)\exp(-x-y)\cX_{[0,1]^2}
$$
given access to samples of the form
$\cF = \{\hat f(k_1,k_2): k_1,k_2 = -512,\ldots, 512\}$. Suppose that we observe only $4.25\%$ the samples in $\cF$, namely, those restricted to the star-shaped domain $\Omega$ shown on the top right of Figure \ref{fig:CS}. Note that the reconstructions obtained by solving (\ref{eq:gsl1_fin}), with $K=256$ and (\ref{eq:cs_fin}) are both given the same input samples, but there is a substantial difference in reconstruction quality due to the samples mismatch introduced by the use of the finite dimensional matrices in (\ref{eq:cs_fin}). Note also that, when subsampling from the first $M$ samples of the lowest frequencies, the computational complexity of computing (\ref{eq:cs_fin}) is $\ord{M\log(M)}$ since the computational complexity of applying $U_{df}$ is $\ord{M\log(M)}$, and the computational complexity of applying $U_{dw}$ is $\ord{M}$. Thus, the computational complexity of computing (\ref{eq:gsl1_fin}) is no worse than the finite dimensional approach of (\ref{eq:cs_fin}).

\appendix

\section{Computation of the Fourier transform of boundary scaling functions}

For the precomputational step in our algorithms, we need to compute the Fourier transform of boundary scaling functions. It is well known that the Fourier transform of scaling functions can be approximated using the low pass filter of the wavelet system. The key ideas on how to do this are recalled in this section.
First observe that  the internal scaling functions satisfy the following equation:
$$
\phi (x) = \sqrt{2} \sum_{n} h_n \phi(2x - n)
$$
which gives rise to the equation
$$
\hat \phi(\xi) =  \prod_{j=1}^\infty m_0(2^{-j} \xi)
$$
where $m_0(\xi) = \frac{1}{\sqrt{2}}\sum_{n} h_n e^{-2\pi i n\xi}$. This provides a recursive way of approximating $\hat \phi$, since $\hat \phi(\xi) \approx  \frac{1}{\sqrt{2\pi}} \prod_{j=1}^N m_0(2^{-j} \xi)$ for large $N$. The same principle can be applied to compute the Fourier transform of the boundary scaling functions, for completeness, we will write out the computations in the case of the left boundary functions.
 
Recall from \cite{cohen1993wavelets} that  the left boundary scaling functions satisfy the following equations for each $k=0,\ldots, a-1$
$$
\phi^0_{k}(x) = \sqrt{2}\left(\sum_{l=0}^{a-1} H_{k,l}^0 \phi_{l}^0(2x) + \sum_{m=a}^{a+2k} h_{k,m}^0 \phi(2x - m)\right),
$$ 
where $\br{H_{k,l}^0}$ and $\br{h_{k,m}^0}$ are the filter coefficients associated with the left  boundary scaling functions (these filters are available from \url{http://www.pacm.princeton.edu/~ingrid/publications/54.txt}).
So, 
$$
\hat \phi^0_{k}(\xi) = \frac{1}{\sqrt{2}}\left(\sum_{l=0}^{a-1} H_{k,l}^0 \hat\phi_{l}^0\left(\frac{\xi}{2}\right) +  \hat\phi\left(\frac{\xi}{2}\right)\sum_{m=a}^{a+2k} h_{k,m}^0 e^{-2\pi i m \xi/2} \right)
$$

Let
$$
U \! = \! \begin{pmatrix}
H_{0,0}^0 &H_{0,1}^0 &\hdots &H_{0,a-1}^0\\
H_{1,0}^0 &H_{1,1}^0 &\hdots &H_{1,a-1}^0\\
\vdots & & &\vdots\\
H_{a-1,0}^0 &H_{a-1,1}^0 &\hdots &H_{a-1,a-1}^0
\end{pmatrix}, \ \ V \!=\! \begin{pmatrix}
h_{0,a}^0 &0  &0 &\hdots  &\hdots &0\\
h_{1,a}^0 &h_{1,a+1}^0 &h_{1,a+2}^0 &0 &\hdots
 &0\\
 \vdots & & &  & &\vdots\\
h_{a-1,a}^0 &h_{a-1,a+1}^0 &\hdots &\hdots &\hdots &h_{a-1,3a-2}^0
\end{pmatrix}
$$
and 
$$
v_1(\xi) = \begin{pmatrix}
\hat \phi^0_{1}(\xi)\\
\hat \phi^0_{2}(\xi)\\
\vdots\\
\hat \phi^0_{a-1}(\xi)
\end{pmatrix}, \quad v_2(\xi) = \hat\phi(x) \begin{pmatrix}
e^{-2\pi a\xi}\\
e^{-2\pi (a+1) \xi}\\
\vdots\\
e^{2\pi i (3a-2)\xi}
\end{pmatrix}
$$
It follows that for each $j\in\bbN$
$$
v_1(\xi) = \frac{1}{\sqrt{2^j}} U^j \left(v_1\left(\frac{\xi}{2^j}\right)\right) + \sum_{l=0}^{j-1} \frac{1}{\sqrt{2^l}} U^l V \left(v_2\left(\frac{\xi}{2^{l+1}}\right)\right). 
$$
We remark also that since $\hat\phi(0) = 1$, $v_2(0)$ is the vector whose entries are all ones and the values $v_1(0)$ can be obtained by solving 
\be{\label{scaling_solve}
v_1(0) = U (v_1(0)) + V( v_2(0)).
}
Note that $\hat \phi_j^0$ is smooth and so, for each $\xi\in\bbR$, $v_1(\xi/2^j) \to v_1(0)$  as $j\to \infty$. So,
$$
\nm{ v_1(\xi) - \left(U^j \left(v_1\left(0\right)\right) + \sum_{l=0}^{j-1} U^l V \left( v_2\left(\frac{\xi}{2^{l+1}}\right)\right)\right)} \to 0
$$
as $j\to \infty$ and for $j$ sufficiently large, we may approximate $v_1(\xi)$ by 
\be{\label{scaling_approx}
v_1(\xi) \approx U^j\left( v_1\left(0\right)\right) + \sum_{l=0}^{j-1} U^l V\left( v_2\left(\frac{\xi}{2^{l+1}}\right)\right).
}
Thus, for a given $\xi\in\bbR$, to approximate $v_1(\xi)$, we first solve (\ref{scaling_solve}) to obtain $v_1(0)$, and then we compute  (\ref{scaling_approx}) for $j$ sufficiently large.
 
A similar approximation can be derived for the Fourier transforms of the scaling functions on the right boundary.

\section*{Acknowledgements}
 This work was supported by the UK Engineering and Physical Sciences
Research Council (EPSRC) grant EP/H023348/1 for the University of Cambridge
Centre for Doctoral Training, the Cambridge Centre for Analysis. The authors would like to thank Ben Adcock and Anders Hansen for their comments.

\addcontentsline{toc}{section}{References}
\bibliographystyle{abbrv}
\bibliography{References}

\end{document}